\documentclass[12pt, english]{article}
\usepackage{latexsym}
\usepackage[activeacute]{babel}
\usepackage{graphicx}
\usepackage{amsmath}
\usepackage{amsthm}
\usepackage{amsopn,amstext,amscd,amsfonts,amssymb}
\theoremstyle{remark}

\theoremstyle{definition}

\begin{document}
\title{A glimpse into the life and times of F.V. Atkinson}


\author{Angelo B. Mingarelli, \footnote{Dedicated to the memory of F.V. Atkinson, teacher and mentor. This research was supported in part by grants from the Office of the Vice-President Research, Carleton University and NSERC Canada.}
\footnote{School of Mathematics and Statistics and Faculty of Graduate Studies and Research,
Carleton University, Ottawa, Ontario, Canada, K1S 5B6,\quad e-mail: {\sf amingare@math.carleton.ca}} }
\date{}

\maketitle

\begin{abstract}
We review the life and scientific work of F.V Atkinson, mathematician. The article consists of a brief overview of his academic and professional life followed by an extensive compilation of his life's work, some unpublished. He is best remembered for his contributions to Number Theory (mean value theorems for the Riemann zeta function on the critical line), his work in operator theory (index theorems for abstract Fredholm operators), and his seminal work in many areas of differential equations (multiparameter spectral theory, Sturm-Liouville theory, asymptotics etc.).
\end{abstract}


{{\begin{verse}{``Even after we understand what they [the magicians]  have done, the process by which they have done it is completely dark ..."\quad
Mark Kac.}\end{verse}}}
\vskip0.15in
\begin{center}CONTENTS\end{center}
\begin{enumerate}
\item Biographical information \hfill \ldots\ldots~\pageref{biog}
\item Atkinson as a teacher \hfill \ldots\ldots~\pageref{teacher}
\item Scientific and scholarly research \hfill \ldots\ldots~\pageref{ssr}
\begin{itemize}
\item The Riemann zeta-function and analytic number theory \hfill \ldots\ldots~\pageref{rh}
\item Asymptotic and oscillation theory of ordinary differential equations \hfill \ldots\ldots~\pageref{aot}
\item Functional analysis and operator theory \hfill \ldots\ldots~\pageref{fa}
\item Functions of a complex variable \hfill \ldots\ldots~\pageref{fcv}
\item Orthogonal polynomials \hfill \ldots\ldots~\pageref{op}
\item Spectral theory of ordinary differential operators \hfill \ldots\ldots~\pageref{std}
\item Partial differential equations \hfill \ldots\ldots~\pageref{pdes}
\item Miscellaneous topics \hfill \ldots\ldots~\pageref{mt}
\end{itemize}
\item Publications list \hfill \ldots\ldots~\pageref{pubs}
\item Miscellaneous articles \hfill \ldots\ldots~\pageref{mas}
\item Impact of the works of F.V. Atkinson \hfill \ldots\ldots~\pageref{five}
\item Miscellaneous references \hfill \ldots\ldots~\pageref{mrefs}
\item Acknowledgments \hfill \ldots\ldots~\pageref{ak}
\end{enumerate}

\section{Biographical information}
\label{biog}

Frederick Valentine Atkinson (known to his close friends as Derick) friend to many, scholar, enthusiastic teacher and gifted researcher closed his eyes for the last time on November 13, 2002, after a long illness, in Toronto, Canada. He was born at Pinner, Middlesex in England on January 25, 1916, the elder son of George Arthur Atkinson and Dorothy Boxer, his second wife, brother to Ann, his only sister four years his junior and half-brother to Victor. It is something of a mystery that Derick only found out that Victor was his brother late in life. From then on he made regular visits to Victor's home in  England and cared for him just as he would for any family member. His father, a journalist and film critic for the {\it Daily Telegraph}, was a promoter during his long and distinguished career of the early works of Alfred Hitchcock. His mother's grandfather was Lord Admiral Boxer, Harbourmaster of the City of Qu\'ebec during the 19th. century.

The earliest document available on Derick's early childhood dates from when Derick, 6 years young, was registered in Preparatory Form at Reddiford School  in Cecil Park, Pinner. A {\it report card} from this school dated April 22, 1922, attests to Derek's (sic) attentiveness in Scripture, his very good {\it Number Work}, his much improved writing and his good memory. The next glimpse of his educational development in the early school years comes in 1925 where another report card now places him at Alleyn's School, finishing third in overall ranking out of a total of 30. His best grades were in English Dictation,  Handwork, Drawing, and Arithmetic (in decreasing order of proficiency). His teachers were very pleasantly surprised at his progress: {\it `` He always seemed to be half asleep in class, yet he must have taken things in,"} wrote one of his teachers. In 1927 he and his friends started up a newsletter of sorts entitled {\it The Forest Rovers}, completely handwritten and with very limited distribution (less than a handful is surmised). Derick (and his sister Ann) contributed to the Sept. 3, 1927,  issue each with a short story entitled {\it The Faithful Dog} (and {\it The Goblin who Told Tales}, respectively). Atkinson closed his story with the simple line ``Derick Atkinson, Aged 11".

It has been said that Atkinson read books on Calculus at age 12. Although this claim is purely apocryphal there is a report from Dulwich College Preparatory School (where Derick was enrolled in the Spring of 1929)  which indicates an excellence at Geometry, Algebra, and Arithmetic, placing at the top of his class in all these subjects.  Presumably this was a local school since at this time his parents lived on Turney Road in Dulwich Village. Although mathematics came easy to the young Derick his parents had different plans for his future; his father wanted him to become a diplomat while his mother would have preferred a clergyman.

Soon after in June, 1929, Derick was accepted into the {\it classics} program of the old and legendary St. Paul's School in West Kensington.  One can guess that these were his most important years as to his future plans. The formative years, 1929-1934, were spent in the same corridors whose walls once reflected the shadows of the poet John Milton, the diarist Samuel Pepys, the mathematicians J. E. Littlewood (of {\it Hardy and Littlewood} fame) and G.N. Watson  among scores of others in a school that was the cradle of  many historical figures in the humanities. It is remarkable that the same John Littlewood would sit as a member of the examining board of Atkinson as the latter would prepare to defend his doctoral dissertation at Oxford years later. By the Autumn Term of 1931 his destiny as a future mathematician was likely established. Indeed, a handwritten report card from St.Paul's praises his mathematical ability to the point where ``... I hope he will keep up interests outside his mathematics and guard against becoming unbalanced." Above the High Master's signature is the final closing simple assessment: `` Extremely promising: He should make a brilliant mathematician". Derick was just over 15 years old at this time. John Bell, the High Master there from 1927-1939, attested to this distinguished period in his early life. In a letter of reference on behalf of Atkinson, Bell enthusiatically sums up his abilities with the statement, ``... [Derick] bore an excellent character and showed first rate ability, especially for Mathematics, in which he won an open scholarship at Queen's College, Oxford". While at St. Paul's he won a book entitled {\it Lectures on Ten British Mathematicians of the Nineteenth Century} by Alexander MacFarlane (1916). It was presented to him by the High Master himself. All this praise bore fruit on the day when in the summer of 1932, George Atkinson was notified that his son Derick had been recommended for election to a Senior Foundation Scholarship.

In November, 1933, at age 17 he was notified by Queen's College, Oxford, that arrangements had been made on his behalf to participate in the Mathematical Scholarship Examination there. In 1934 he enrolled at The Queen's College. During his stay at Queen's he once shared a room with a close future friend, the late mathematician J.L.B. Cooper. He was a Tabardar at Queen's, so called because the formal gown worn by all had tabard sleeves, {\em i.e.,} loose sleeves, terminating a little below the elbow in a point. The year 1935 saw him achieve First Class Honours in the Mathematical {\it Moderations}  (the first part of the degree, - it was necessary to pass Moderations at the end of the first year in order to pass on to the Honours exam). For this achievement he received a note of congratulations from the Provost of Queen's, Canon B. H. Streeter. As a further consequence of this accomplishment he was invited to a special {\it Mods. Dinner} on $7^{th}$ March, 1936. The letter of invitation came from one silently influential figure in Derick's years at Oxford,  C. H. Thompson (Fellow of Queen's College, 1897-1948) bearing the date  13.ii.36 on Queen's College letterhead.  The reply to this note by Derick reveals some interesting characteristics of the young brilliant mathematician. The note reads: {\it  Dear Atkinson: Will you give Mr. Wright and me the pleasure of your company at the Mods. Dinner in the Upper Common Room at 7:45 p.m. on Saturday, 7th March. Yours Sincerely, C.H. Thompson}.  Derick's reply on the verso was: {\it Dear Tompey, Just to please thee, I mussn't talk to that theer Wright. 'Es a bugger he be. Yours affectionately, F.V. Atkinson}. It is curious that the original note somehow made its way {\it back} into Derick's hands ... It is possible that the {\it Mr. Wright} referred to here was in fact, T. E. Wright, Fellow of Queen's College (1929-1948), - he may have been one of Derick's tutors at Queen's.

On the March 24, 1936, Derick was notified by the Provost of the College that he had obtained an Honourable Mention of sorts, lit. {\it proxima accessit}, for the Junior Mathematical Exhibition Examination. This was followed by a note from the Governing Body of the College attesting to this and to the ``{\it award of a prize of books to the value of  $5\pounds$}".  He was also awarded a First Class Final Honour School of Mathematics in 1937, the same year he obtained his B.A. in Mathematics with {\it First Class Honours} from Queen's. The name of his long time friend and colleague, J. L. B. Copper, also appears on the Honour Roll for 1937. Very little is known of the period 1937-1939, the years in which he spent working closely with his mentor, E.C. Titchmarsh, who held the Savillian Chair at the time (once held by G. H. Hardy, his own supervisor). There is, however, a typed copy of a letter of reference from Titchmarsh in Derick's original ``job search kit".

\begin{quote} I have much pleasure in recommending Mr. F.V. Atkinson for a position as mathematical lecturer. He has worked with me as an advanced student for nearly two years, and he has done some interesting work in the analytical theory of numbers. I feel confident that he will obtain the D. Phil. degree in due course. He has worked very hard and shown originality and resource in tackling some very difficult problems. He would be personally an agreeable colleague to work with.

\flushright (Signed) E.C. Titchmarsh

\flushleft

Feb. 28, 1939.

\end{quote}

Derick also had a copy of another letter of reference, one by U. S. Haslam-Jones, a Fellow and Tutor of Queen's College dated March 4, 1939. This fascinating letter, part of which is reproduced below,  gives a rare astute glimpse into a trait of Derick's that was almost indiscernible  ...

\begin{quote} ... He is primarily an analyst, but his work in Applied Mathematics for me was always sound. In his undergraduate days he suffered, to some extent, from a fault which is almost an undeveloped virtue - namely, an eagerness to get on with his subject which led him to leave a problem as soon as he had seen how it should be solved, and to dislike the drudgery of working on it until the very end. His years of research seem to have cured this fault, and he is now far more methodical...

\flushright (Signed) U. S. Haslam-Jones

\flushright Fellow and Tutor of Queen's College

\end{quote}

Derick's tutor while at Queen's was C.H. Thompson (``Tompey" mentioned above). It is a curious twist of fate that the only other senior mathematician at the College during the period that Derick was there was A.E.H. Love (Augustus Love) who held the Sedleian Chair of Natural Philosophy. Recall that Love is that one mathematician who in effect  inspired G.H. Hardy (cf., \cite{hardy}) to read Camille Jordan's {\it Cours d'Analyse}, a book that would eventually lead to Hardy's own monograph ({\it A Course of Pure Mathematics}) thereby changing forever the face of pure mathematics in the UK. It is more than likely that Love and Atkinson knew each other, in a sense Augustus Love was his mathematical great-grandfather! While at Queen's College he took his D. Phil in $1939$ under the supervision of Titchmarsh, then the {\it Savillian Professor of Geometry}.  His dissertation consisted in finding asymptotic formulae for the average value of the square of the Riemann zeta function on the critical line, \cite{thesis}.  An abstract of his D.Phil. dissertation may be found in \cite{abstract} and, at this time, only a draft of the final dissertation exists (totalling 68 pages, plus one page of references).  He was fond of recounting to interested students that his final Examining Board consisted of  G.H. Hardy, J.E. Littlewood and E.C. Titchmarsh! While at Oxford it has been confirmed \cite{archives}, that he was secretary of the Chinese Student Society and also a member of the Indian Student Society. Although not well known, he was fluent in many languages including English, Latin, Ancient Greek, Urdu, German, Hungarian, Russian with some proficiency in Spanish, Italian and French. His fluency in Russian was all the more remarkable given that he never (officially) followed a Russian course at Oxford, \cite{archives}!

Upon completion of his D.Phil. program he was appointed Senior Demy (a Junior Research Fellow of sorts) at Magdalen College, Oxford, during the academic year $1939-1940$ but this appointment was cut short by World War II. In a letter to the President of Magdalen College dated October 5, 1940, Derick felt obliged to resign his Senior Demyship. The President clearly regretted the sudden departure; in a reply dated October 8, 1940, he says: `` ...I realize, of course, that the claims of your new work must override at the present time all other considerations. I am glad that the College has been able to assist you at what would otherwise have been an important moment in your work and research." It is likely that the inference to the {\it claims of your new work} was referring to Derick's acceptance of a commission into the Intelligence Corps that would become effective December 12, 1940.

Sometime during the Autumn of 1940, every Tuesday at their home in Oxford, Derick began visiting a family originally from Slovakia, Ladislav and Agnes Haas, both physicians. Ladislav was a psychiatrist trained in the freudian tradition who was redoing his MD in England and Agnes was a pediatrician. They had quietly left Slovakia for England since he was being sought by the SS for his political views. Derick wanted to practice his Hungarian and Russian with the Haas family and they, in return, would practice their English with him since their medical degrees would require re-certification in order for them to practice in England.  Ladislav and Agnes had two daughters, Dusja and \v{S}ura. Dusja, the youngest of the two,  was born in Kuybishev, in the former Soviet Union. He had met Dusja, who at the time spoke almost no English, through her father. As an aside, Dusja's mother's first cousin married Salo Flohr, a  famous Chess Grandmaster of the period. The once World Chess Champion, Michail Botvinnik, once compared Flohr to Napoleon and said: ``In the 1930s they all trembled for Flohr". In the early days Dusja and Derick spoke Hungarian when together, a language in which Derick was fluent. After a very brief engagement he married Dusja Haas in the Registry Office in Slough, England (between London and Oxford). Such was Derick's modesty that it was only on this day that Dusja discovered that her husband had a doctorate! Dusja and Derick had three children, Stephen, Vivienne and Leslie, the latter two having joined the ranks of academe as professional psychologists following, to some extent,  the path of their maternal grandparents.


During the war years he worked in the Government Code and Cypher School at Bletchley Park. However, in 1946 the UK Government ordered the destruction of most records held at Bletchley Park including all personnel records. Nevertheless, a small amount of information on Derick can be gathered from documents subsequently released by the Public Records Office at Bletchley Park. F.V. Atkinson is mentioned in only one document in existence and it is a list of people working at Bletchley Park on December 2, 1941. On that day he was stationed at Wavendon, a mere 4 miles from Bletchley. Wavendon was first established to act as a {\it duplicate} site in the event of air-raid damage to Bletchley Park, \cite{bp}. The first Bombe outstation known as OSW (Out Station Wavendon) was established in early 1941, \cite{bp}. No record of the duration of his employment at Wavendon exists, however. Although little is known of this period of his life because of the Official Secrets Act and because of the official destruction of all personnel files, it is likely that he knew the famous computer scientist Alan M. Turing, who worked in the same location at the same time. He was promoted to Captain on January 1, 1943 and achieved the rank of Major in the British Intelligence Corps sometime in 1945.

Derick's noble, generous and unselfish nature can be gathered through a surviving document arising as a result of the rationing of sugar and tea during the war years in England. The following excerpt of a letter dated Feb. 22, 1944, from Titchmarsh to Atkinson gives an impression of the academic environment at Oxford during the war years:

\begin{quote}Dear Atkinson,

Thank you very much for the parcel which arrived safely a day or two ago. It was extremely kind of you to think of us after all this time. We are already enjoying home-made marmalade made with your sugar.

I hope you get some interesting work to do in India, though I never heard just what your line was. Everything here is very dead. There are about six women candidates for maths finals this time and no men. I am doing some very elementary teaching - it is years since anyone came to my advanced class ... I hope we shall see you again before so very long, but it must be rather difficult to get back from India. ....

(signed) Yours sincerely,

E.C. Titchmarsh

\end{quote}

This letter must have been written with the full knowledge of the authorities involved. Indeed, Atkinson was in India at this time (for a total of about three years) helping to break the Japanese code for British Intelligence. Notwithstanding the fact that Derick wrote very few mathematical papers that we know of during this period of War Service that is, from 1940-1946, he was formally offered an appointment on December 28, 1945, (later declined) as an Assistant Professor in the Department of Mathematics at Marischal College, in the University of Aberdeen. This appears to be his first formal offer of an academic appointment outside of the Oxford Colleges. Existing documents indicate that on the January 2, 1946, Derick was definitely part of the Intelligence Pool in the Intelligence Corps Depot at Rotherham in Yorks where he had been posted to the position of {\it Military Technical Officer Class II} in the Army Operational Research Group. He was living in Oxford at the time.

His release from active military duty in the Intelligence Corps occured on October 3, 1946, with a letter of thanks from one Eric B. B. Speed, on behalf of the Army Council in the War Office:

\begin{quote}Sir:

Now that the time has come for your release from active military duty, I am commanded by the Army Council to express to you their thanks for the valuable services which you have rendered in the service of your country at a time of grave national emergency.

At the end of the emergency you will relinquish your commission, and at that time a notification will appear in the London Gazette (Supplement) ...

\end{quote}

Still, in a confidential letter dated October 22, 1946, signed by C. G. Hobden, Assistant Adjutant-General, Derick was asked to be part of an {\it active} list of  ex-intelligence officers in the Post-War active Army. There is no documentary evidence at this time indicating whether or not Derick accepted this commission although it is suspected that he may have since he wore his military uniform while lecturing at Oxford during this period.

The details leading to his being awarded an M.A. from Oxford sometime in $1946$ are obscure and the subject of  speculation. On the July 12, 1946, he was elected to a Lectureship in Christ Church, Oxford, his first (accepted) formal academic appointment that would terminate in 1948. A note dated December 9, 1947, from the University Registry, Oxford, requests  that Derick lecture twice per week on Special Functions in {\it Hilary Term} and once per week on Elliptic Functions in {\it Trinity Term}. Compensation seemed to depend upon the number of lectures given (a practice that dates to medieval times).

In a note dated June 16, 1948, Derick's Lectureship in Christ Church was extended for a further two-year period. However, scarcely a week later, in a formal letter dated June 22, and originating from the Office of the Under Secretary of State, he was offered a (Full) Professorship in Mathematics at University College, Ibadan, in Nigeria. It was an offer that would be difficult to refuse. For instance, the terms of the offer included, in part, a yearly  {\it expatriation allowance}, in addition to his nominal salary, because he was not normally domiciled in West Africa. He would be provided with rental accomodation that would not exceed $2$\% of his basic salary, and he would be entitled to free medical services and hospital treatment at a reduced rate. In addition, it was an appointment with tenure and he would have the ability of traveling first-class while {\it on duty}.  His resignation as a Lecturer at Christ Church became effective on the September 29, 1948.

His duties while in Ibadan (albeit onerous ones) can be found in his letter of appointment:

\begin{quote}

Your duties will initially be to develop your department as a centre of teaching and research. This will include discussion with the appropriate members of the University of London of the syllabus of your subject within the terms of the special relationship with the University of London; assistance to the architect in the design and building of the permanent accomodation for your department; and, in consultation with those in charge of the library, the compiling of book-lists and the taking of other steps to arrange for the provision of library facilities in your subject. You will be a member of the Academic Board of the College ...

\end{quote}

Besides being Head of the newly created Department of Mathematics at University College, Ibadan, Derick's research shifted (for no apparent reasons) from Number Theory to Functional Analysis, Operator Theory and Differential Equations. It was here that he wrote some of the most enduring papers of his youth, \cite{1949-4}, \cite{1951-1}, \cite{1952-1}, and \cite{1953} that helped to bridge the gap between differential equations and operator theory.

He was appointed Dean of the Faculty of Arts in Ibadan during the academic year 1953-1954 and in addition, he taught Statistics and Physics courses as part of his regular teaching load. His work not wholly devoted to the academic environment during this academic year we find his nomination to the Fellowship of SIRA, the Argentine International Broadcasting Service, for spreading the news and relating events pertaining to life in Argentina. His interest in Statistics appears to be more than simply pedagogical. Through the years he contributed to both the theory and applications of statistics. For instance, he obtained a citation in a paper on anaesthesia \cite{anes}, and he actually co-authored a paper in statistics, \cite{1964-1}. In his spare time he could be found occasionally playing the accordeon in a local club for free beers! Although life seemed blissful, the heavy teaching load during his stay in Ibadan (more than 12 contact hours per week) compounded with administrative duties of every kind eventually took its toll on Derick and thus led him to look elsewhere.

It is difficult to determine the sequence of events following his departure from Ibadan (Spring, 1955) but preceding his arrival in Canberra, Australia (Autumn, 1956). What is certain is that he was made a formal offer (Summer, 1955) of membership in the Institute for Advanced Study at Princeton for the academic year 1955-1956 subsidized in part by a contract from the Office of Ordnance Research. Some say he was invited there by Leopold Infeld, Albert Einstein's collaborator in his later years. In a letter dated August 27, 1955, the Institute had actually taken the initiative to reserve a room for him in anticipation of his arrival there in a building that had been housing other Institute members. His failure to take up this position appears to have been due to his inability to secure a visa for the US. The Atkinsons were living in Derick's parents' house in Sussex, England during this academic year and Derick visited Oxford University regularly. He earned his basic living to support his young family by translating mathematical papers from Russian into English (see e.g., \cite{bautin})

He moved to Canberra University College (now part of Australian National University, or ANU as it is referred to nowadays) in $1956$ to become Head of its Department of Mathematics. The Atkinsons were met in Canberra by the Principal of the College, Dr. Joe Burton who eventually became a dear friend of the family. While there he made frequent administrative visits to the University of Sydney to review examinations and course offerings. In response to a letter from Prof. M. L. Oliphant \cite{oliphant}, of Canberra requesting Atkinson's advice on the establishment of a new Department of Mathematics at  ANU, Derick had the following words of wisdom:

\begin{quote}\quad\quad ... What you say about mathematicians being ignorant of each other's specialties is all too true; I think, however, that this is less due to the actual diversity of the specialties than to the preoccupations of a teaching post. Some of the really great mathematicians have been able to make important contributions to several fields; the late H. Weyl is an example. The gap between geometry and the theory of numbers seems to have been narrowing during the past two decades. For reasons of this kind I would consider it impracticable to tie down the posts to particular subdivisions of mathematics; the most that could be done in this direction would, I think, be the creation of posts in pure and applied mathematics, though even this classification is apt to break down, and I am not advocating it.

\flushright{Letter to M. L. Oliphant, $13^{th}$ February 1958.}

\end{quote}

Sometime in November $1959$, he wrote a letter to the {\it Saturday Evening Post} that portrayed his passionate interest in mathematics, one that aroused some  interest from its readership. It was in Canberra that he met his first Ph.D. student, Claude Billigheimer, who eventually followed him to Canada when Atkinson left Canberra for the University of Toronto in $1960$. His departure from Canberra College was due in part to an excessive teaching load and other onerous administrative duties.

We outline a few details of his departure from Canberra for his final home in Canada. It appears that Derick had applied for positions at both Queens University in Kingston, Ontario, and the University of Alberta in Edmonton, Alberta. At a meeting in Queens Dr. Jeff Watson (then at the University of Toronto) met with an administrator from the Department of Mathematics in Kingston who had promptly informed him that one, F.V. Atkinson from Australia, was just made an offer of an Associate Professorship. Before the conversation could end Watson left immediately by car for Toronto (more than a two hour drive away) and discussed this matter with then Chair, D. B. De Lury. Some faculty members there already had a connection with Canberra University College and knew of Derick's reputation. Soon after there followed a cable to Derick Atkinson in Australia offering him a position as Full Professor at the University of Toronto, a position that was accepted, effective  July 1, 1960. The Atkinsons arrived in Canada by ocean liner on October 1, 1960, via Montreal. This was followed by a short train journey to Toronto. Since the semester had already initiated at the university it became the responsibility of Dr. P.G. (Tim) Rooney to teach in Atkinson's place for part of  that term.  This same year he took up the post of {\it Associate Editor} of the Journal of Mathematical Analysis and its Applications, a position he would keep for nine more years.

He spent the summer of 1961 at the Functional Analysis Summer Institute at Stanford University, the summer of 1962 at the Canadian Mathematical Congress Summer Research Institute at Queen's University, in Kingston, Ontario, where he wrote \cite{1963-1}, and the summer of 1963 as a Visiting Professor at the Mathematics Research Center at the University of Wisconsin, in Madison. The year 1964 saw the publication of his next {\it opera omnia} \cite{1964}, a seminal work composed and written completely in Canada, a monograph that would be translated quickly into Russian (1968) with notes by famous Russian mathematicians of the day. Throughout his career, Derick was always in much demand from his colleagues around the world; it is said that Richard Bellman tried to attract him to the University of California at Los Angeles, and W.C. Royster, then Chair of the Department of Mathematics and Astronomy at the University of Kentucky tried to persuade him to come there (an informal solicitation does exist). This last offer was declined gently although Derick did spend the academic year 1964-1965 in Kentucky.

Permanently positioned in Toronto he served as Secretary of the Graduate Department in 1963 for a few years. He accepted a position as Associate Editor of {\it Mathematical Reviews} in 1965. He graduated his first two Ph.D. students, Claude Billigheimer (formerly of McMaster University in Canada) and Melvin Faierman (formerly of the University of the Witwatersrand, South Africa) in 1966. He was elected Fellow of the Royal Society of Canada in 1967. After being a Visiting Lecturer at the University of Freiburg in 1969 he took up a U. K. Science Research Council Visiting Fellowship at the Universities of Dundee and Sussex in 1970. His third doctoral student, Stephen Ma graduates in 1972. In 1973 he was made a British Council Lecturer to U. K. universities. Derick was made Consulting Editor of the {\it Proceedings A} of the Royal Society of Edinburgh in 1974, the same year that he graduated his fourth doctoral student K. S. Ong. He held this editorship for ten years.

In 1975 he was elected Honorary Fellow of the Royal Society of Edinburgh. As an aside, this honor deserves more than just a cursory reading. In 1991, there were only 64 such Honorary Fellows; the list included the mathematicians  Michael F. Atiyah, Mary L. Cartwright (another student of Hardy/Titchmarsh), Gaetano Fichera, Jack K. Hale, Paul R. Halmos, Jaroslav Kurzweil, Saunders MacLane and Olga Oleinik. At that time, this honorary membership list included 13 Nobel Laureates and consisted mostly of fellows of the Royal Society of London.

Derick was elected Chair of the University of Toronto's Department of Mathematics in $1975$ a position he held for six more years. He was a member of the Advisory Council of the Canadian Society for the Weizmann Institute of Science (ca. 1976). In 1977 Atkinson was awarded the {\it Makdougall-Brisbane Prize for 1974-1976} of the Royal Society of Edinburgh for his paper \cite{1974-2}. Other notable recipients of this Prize include the Nobel Prize-winning physicist, Max Born, and the mathematician  J.H.M. Wedderburn. He graduated his fifth doctoral student M.W.Wong in 1978 and his last such student, the author, in 1979.

He was appointed Visiting Professor and Consultant (with NSF support) at the Pennsylvania State University in University Park, in 1980. In 1981 Derick was appointed to the Grant Selection Committee for Pure and Applied Mathematics at N.S.E.R.C. (Natural Sciences and Engineering Research Council of Canada) for a period of three years. In $1982$, the year in which he was made (U.K.) Science Research Council Fellow at the University of Dundee in Scotland, he retired as Professor Emeritus from the University of Toronto. After this event there is a clear discontinuity in his research output and the next 10 years of his mathematical life are marked by tremendous bursts of scientific research activity: Almost one-half of his lifetime total of published papers appeared during this ten year period. During the summer of 1983 he visited Argonne National Laboratory in Argonne, Illinois, while in 1984 he was Chair of the Unit Review Committee at the University of Alberta and a Visitor in the Department of Mathematics there. He took on Visiting Professorships at the University of Alabama in Birmingham in 1981, 1983, and 1984. He visited the University of Leiden in the Netherlands frequently and was S.E.R.C. (Sciences and Engineering Research Council) Visiting Research Fellow at the University of Dundee, Scotland, in 1985. He also served the Canadian Mathematical Society for more than $40$ years. Among the many duties he undertook on its behalf we cite his tenure of the office of  its $29^{th}$ President, $1989-1991$ and his former membership on the Board of Directors, $1967-1969$.

We outline his post-retirement chapter (1985-1992) briefly: He spent much of his time in visits to Argonne National Laboratory (Hans Kaper and others), the University of Leiden (L. Peletier), and the University of Regensburg (Germany, R. Mennicken and others). He also visited the Technical University of Vienna (H. Langer), and the universities of Winnipeg and Calgary, in Canada as well as the University of Maryland in the USA. He was a keynote speaker at an {\it International Conference on Boundary Value Problems for Ordinary Differential Equations} at Ohio University in Athens, OH and at ICIAM (1991) in Washington DC. In 1991 he attended an {\it International Conference on Nonlinear Waves} honoring his $75^{th}$ birthday at Carleton University in Ottawa. During the stated period he held one of the largest NSERC Research Grants in all of Canadian mathematics. Among the many additional honors bestowed upon him we cite the Barrett Lectureship at the University of Tennessee (also awarded to Einar Hille) and the {\it Von Humboldt Prize} by the government of Germany (1992). That same year, in August, he suffered a massive stroke followed by other complications upon his return from Europe that disabled him completely until his death in $2002$.

\section{Atkinson as a teacher}
\label{teacher}

By way of personal recollections of Derick and his influence on many a student, I can recall the days when, as a junior undergraduate in a small Jesuit college in Montreal, I had Mel Faierman as a teacher in a course on Advanced Calculus (recall that Mel was Derick's second doctoral student). Mel was a definite influence in shaping my future aspirations as a mathematician. I remember him saying one day: `` If you {\it really} want to learn mathematics you should go to Toronto...". Of course, he meant the University of Toronto which, at that time, had more than one hundred mathematics graduate students and an outstanding faculty complement including two FRS (Fellows of the Royal Society of London), namely Heilbronn and Coxeter. Mel continued ... ``...and if you {\it really} want to learn differential equations then you should study with Atkinson. " Of course, I didn't know who Derick was at this time but Mel had mentioned that he had studied with Titchmarsh at Oxford in the thirties, so I thought this was cool since I had consulted Titchmarsh's classic text in the theory of functions on a few occasions on the advice of Eric O'Connor, S.J., the former Secretary of the Canadian Mathematical Congress (now Society). On the other hand, I had already been in contact with Coxeter on account of some new results in classical geometry that I had discovered along with a childhood friend, Mario D'Angelo. So I was really thinking about working with Coxeter later on.

Then came the first term of my first year at the University of Toronto. I had courses with the world famous geometer, H. S. M. Coxeter, F. V. Atkinson, and later on with the number theorist, Hans Heilbronn (who had proved the legendary Gauss Conjecture), among others. It was a very difficult term adapting to this new city that was Toronto, the adjustment of a new life as a graduate student, my first time away from home and living alone, etc. In addition, I was learning the highest level of mathematics ever. It was exciting but it was all so difficult to absorb it all. Atkinson's first lecture was most impressive. He walked into a class of about 25 (mostly silent students) wearing a greyish patterned jacket and navy blue pants, a white shirt and tie, his longish silvery hair and penetrating blue eyes making for one who inspired trust along with a commanding presence.

He smiled at everyone and nodded and began writing a list of about ten textbooks that would be recommended reading material for this course, known locally as MATH 435. It was {\it his} course and everyone knew it; this was also Derick's favorite graduate course. The book list included the classic texts by Coddington and Levinson, Hille's lectures on the subject, Hartman's book on the subject, Hale's book on the subject and others. I must admit that consulting these books was a scary experience back then, but it was worth it. He began lecturing on the Banach fixed point theorem (although he called it the fixed point theorem of {\it Banach-Cacciopoli}). He introduced the notion of a Banach space with plenty of examples and everything one needed to understand the statement of this key theorem. Somehow, in this 50 minute lecture, he managed to actually {\it prove} the basic existence and uniqueness theorem for first order differential equations in Euclidean spaces. It was all there, it was clear, and it was complete. I realized immediately that he was a great teacher. I worked very hard at understanding differential equations under his tutelage and it became fun.

One day, sometime during the second half of the first term, at the end of his class, he walked uncharacteristically up the aisle towards me ... I noticed he was carrying a bunch of sheets under his arm and he was smiling as he walked (Derick almost always smiled, it was one of his most endearing traits). As I tried to put my things together without acting clumsy (since I was nervous) he gently placed these few stapled sheets on my desk and walked off. That was it! Well, I thought he had dropped them there and so I gathered my things and ran off towards him as he walked to his office in Sid-Smith (this is what we called the building, actually called Sidney Smith, housing the mathematics department). When I caught up with him, I suggested that he may have dropped these sheets when he walked by ... but he smiled and said ``No, I want you to look at them and tell me what you think...". So I paused, looked at them, thanked him and walked away. I saw a title, no authorship, an introduction, a few definitions and a key theorem. As I looked on I realized the proof was incomplete. Basically, it was a three page paper (I call them {\it sort-of-papers}) with the main proof ending abruptly in the middle of page three. The argument just stopped! As I thought about this incident that evening I suddenly realized that Derick was probably asking me to {\it complete the proof}. I thought this was his way of introducing me to research. I had always wanted to discover new things so this was my chance. It was very exciting! I worked long and hard on this ``proof", trying to complete it. It took me a few weeks but eventually I thought I had it done and even managed to extend the result in some direction. I showed Derick my completed argument and he was visibly pleased but said little else. He knew how to do it, I'm sure, but I guess he wanted to see if I could and to what extent I could.

A few weeks later, closer to the end of term, the same kind of incident occured anew, this time with a different kind of paper and in a totally different area. It was the same kind of setup. A title, no author, an introduction, definitions, something of a proof and then a sudden stop in the middle of an argument. Of course, I knew what to do with the material this time (so I didn't chase after him) and it would take me another few weeks of hard work to complete the proof. By the time his course ended I had actually partaken in the act of original research, - It may not have been entirely original since Derick certainly knew how to do it, but it was new to me and it was really fun doing it.

Now I had to decide who to work with after I completed my Masters degree in 1975. This was difficult since, on the one hand, Donald Coxeter was an outstanding geometer, teacher and humanitarian and I just loved geometry. In addition, Coxeter had helped me to cope with my stressful first term by listening on many an occasion and suggesting things that really helped. Then there was Atkinson who essentially stood alone, a giant among the faculty, one who had showed me how exciting research in analysis and differential equations could be. I thought I would approach each one in turn and ask them what I could work on as part of my doctoral dissertation. Coxeter suggested a problem having to do with sphere packing in eight dimensions while Derick suggested my working in spectral theory of the Sturm-Liouville equation. But then, Coxeter had many doctoral students in various stages of their thesis preparation while Derick had none that year. In looking back I guess I couldn't resist the fact that Derick alone had really fostered in me the excitement for research and I just couldn't see myself working in any other area. So I decided to work with him.

The period 1975-1977 was full of similar incidents with these sort-of-papers. Over that two-year period I must have read four or five papers by Atkinson and I was always asked to submit my comments to him. In few cases the papers were complete manuscripts (still without authorship) but complete with references and, this time, complete proofs. I made suggestions, found typos, fixed a few constants here and there. For example, I actually remember reading the paper \cite{1978-1} before it was ever submitted and before anyone else ever knew of its existence! Those results were really fascinating and I felt very honored to have read these drafts.

This was the technique Derick used in order to teach me to think about differential equations. I eventually understood many of the articles Derick had written during this period. My own dissertation was in yet another area of research; actually it was based on that {\it second} sort-of-paper Derick had given me in class, back in 1974. By the way, I found out eventually that the first sort-of-paper was actually based on his famous publication \cite{1974-2} (but I was given the case $n=2$). As we mentioned above the appearance of that contribution was instrumental in his being awarded the Makdougall-Brisbane Prize of the Royal Society of Edinburgh.

One could relate many wonderful stories about Derick's combined kindness and genius but there is this one story that needs to be retold for posterity. It occured in 1976 when Derick gave lectures in a {\it Seminar in Spectral Theory} that he ran at the university. One day he was reading this paper that had appeared in the {\it Annalen}, explaining everything with a clarity that made the paper actually look simple, when he came upon this conjecture by the author, something about the order of an entire function being precisely 1. He paused, looked down at the page and gently moved the reprint up and down a few times as he gazed its contents in complete silence. Then he walked over to the extreme right-side of the black board where he began to write in very, very tiny symbols something that the rest of us could not follow, and all this in a most serene and pensive mood. Of course we stopped writing and really didn't know what to do with this most uncharacteristic event in our days with Atkinson. Still, with his back turned away from us he wrote and thought, erased and thought, paused and thought and then stopped. All this in maybe four or five minutes at the very most. He then walked back to the center of the board looked at us, smiled and said: ``Yes, he's right" and he continued to read the paper as if nothing had happened. What happened? We all believed that he somehow had actually proved the author's conjecture out there on the very edge of the blackboard and he {\it knew} it. He may have thought it was done now so we should go on with the paper ... We should now recall the prophetic words of his tutor Haslam-Jones (cited above).

Once I became Derick's student I was invited to the many social events at the Atkinsons; their dinner parties were legendary. Derick's wife Dusja always managed to feed and help host these parties almost single-handedly. She was a one-person catering service that could provide for as many as 20 or 30 people at a time. Derick and Dusja became the closest ties that I had to family in Toronto. They were always very kind, understanding, caring, and warm and Dusja was a symbol of tenacity and devotion, strength and hope. She cared completely for Derick until his dying day.

The mathematician Mark Kac once said ``There are two kinds of geniuses,- An {\it ordinary genius} is a fellow that you and I would be just as good as, if we were only many times better." The other kind of genius Kac called {\it magicians}, ``Even after we understand what they have done, the process by which they have done it is completely dark". As an example, Kac placed the famous physicist Richard Feynman in this category of magicians (or wizards). In my view,  Derick Atkinson was a wizard as well. Among his many past-times he enjoyed carpentry, world affairs, and gardening. Atkinson's mathematical genealogy may be viewed on the world wide web as part of the current Mathematics Genealogy Project, \cite{genea}. The gap in our lives once filled by Derick will be impossible to fill, - He was a one-person, one-time phenomenon.

\section{Scientific and scholarly research}\label{ssr}

With few exceptions we note that one of the themes underlying the published research of Atkinson during his 50 years or so as a mathematician was his use of {\it asymptotics} whether they be spectral or functional. He had a mastery of the subject that few possessed. His life's work is dispersed (in his own words) among the broad mathematical fields that include the Riemann Zeta-function and Analytic Number Theory, Functional Analysis, Topics in Functions of a Complex Variable, Orthogonal and Extremal Polynomials, Asymptotic Theory of Differential Equations, and miscellaneous topics in Pure and Applied Differential Equations (both ordinary and partial). The result of recent research indicates that Derick wrote many drafts of papers in traditional applied mathematics while he was in Ibadan, but most of these remained either incomplete or unpublished.

Concerns over the expected length of this paper have resulted in our giving the briefest of summaries of each of  Atkinson's articles in the hope that the interested reader will refer to the actual paper(s) for further details. For additional information on the impact of Atkinson's work in various fields of science the reader is encouraged to refer to Section~\ref{five}, below.

\subsection*{The Riemann zeta-function and analytic number theory}\label{rh}\vskip0.25in

His earliest results here dealt with the finding of asymptotic estimates for the average value of the square of the Riemann zeta-function on the critical line, see his dissertation \cite{thesis}, and the published papers \cite{1939}, \cite{1941}, and \cite{1948-2}. We outline his work in the sequel; references in boldface are followed by brief reviews of that particular reference:

Regarding {\bf \cite{1939}:}\,\, For recent work based on Derick's first published paper \cite{1939}, see \cite{bala1}, \cite{jutila1}, \cite{weinstein1}, \cite{gotze} and the references therein. This research inititiated from Derick's D. Phil. dissertation and culminated in his classic papers \cite{1949-1} and \cite{1950} (more on these below). In his paper {\bf \cite{1939-1}}  there is a representation of two infinite series involving the partition function one of them involving Bessel functions. Recall that interest in the classical partition function was widespread especially after Ramanujan's work on it culminating in Rademacher's asymptotic estimate for $p(n)$ in 1937.  He thanks his teacher Titchmarsh for assistance in simplifying some integrals in \cite{1939-1} and for proving the first two lemmas there (see p. 176 in \cite{1939-1}). For an extensive survey of the partition function that includes Atkinson's contribution, see \cite{gupta}. The article {\bf \cite{1941}}, related to \cite{1939}, has had a slightly wider readership than its counterpart. Some work on improving the error-terms present in \cite{1941} is very recent, {\em e.g.,} see \cite{ivic1} and \cite{hb1}.

The paper {\bf \cite{1941-1}} exhibits an asymptotic form of a partial sum of the function $d_k(n)$, defined as the number of ways of expressing the integer $n$ as a product of $k$ factors. Recent work on this problem may be found in \cite{hb2} and \cite{jwp}.  Judging from its date of submission (June 1941), this appears to be the only article that Derick submitted (maybe even wrote?) while he was actively serving in the Intelligence Corps of the British Army. We note that he gives his affiliation here as ``Oxford" rather than the more secretive {\it Bletchley Park} or {\it Wavendon} clearly for  reasons of national security. Article {\bf \cite{1948}} deals with the Abel summation of some Dirichlet series; some progress here has been made in \cite{chak1} and \cite{bbc} among others. On {\bf \cite{1948-1}}: This is the first of Derick's papers in a language other than English. It contains an approximation of $\log \Gamma (z)$  away from its poles modulo a logarithmic function and a finite Laurent sum containing  Bernoulli numbers. After these few diversions into separate though related areas of mathematics, Derick returns to studies involving the mean-value of the Riemann zeta-function on the critical line with {\bf \cite{1948-2}}. In this paper he finds a lower bound for the number of zeros of the Riemann zeta-function on the critical line for all sufficiently large values of the imaginary part of the argument, in the same spirit as an older result of Titchmarsh dating to 1934.

The classic paper {\bf \cite{1949-1}} reveals Derick's wizardry at finding asymptotic estimates. It is a deep and much quoted paper in the field of number theory dealing with the classical Riemann hypothesis (see Section~\ref{five} and the references therein for more details). In this paper Atkinson improves on remainder estimates of Hardy-Littlewood, and  Littlewood-Ingham-Titchmarsh for the mean value of the square of the Riemann zeta-function on the critical line using techniques akin to those found in his paper \cite{1939}. {\bf \cite{1949-2}} is the first of Derick's papers involving a co-author, Lord Cherwell \cite{lord}, and the first of his post-war articles. It was actually submitted after his release from the Intelligence Corps. In its postscript Derick remarks that the first two theorems therein are due mainly to Lord Cherwell while the remaining three are essentially his. This is likely the only instance in Derick's career where there is a specified division of credit in his joint papers. The paper consists of sufficient conditions for the validity of certain reciprocal relations between arithmetical functions. The authors show that the results include many of the standard mean value and distribution formulae in number theory.  How Derick came to know of  Lord Cherwell (later Viscount Cherwell)  let alone write a mathematical paper with him is unknown. All that we know in this vein is that Lord Cherwell had a Chair in Experimental Philosophy (Physics) at Oxford (during Derick's years there) until he retired from academic life. For more information on the life and times of Lord Cherwell as well as an account of his interests in Number Theory see [\cite{earl}, pp. 99-100]. Finally, {\bf \cite{1950}} exhibits an interesting relation between the Gamma function, the Riemann zeta-function and various Bessel functions. Atkinson's interests shifted mainly towards differential equations around $1950$. Nonetheless, it seems evident that by $1951$, the year of publication of Titchmarsh's classic text on the Riemann Zeta-function  \cite{tit1}, Atkinson was already an authority on the subject as Titchmarsh's acknowledgment in \cite{tit1} attests. There is even an existing copy of his thorough review of a classic paper by A. Selberg \cite{selberg} in number theory.

After an absence of more than a decade, Atkinson returns to number theory via continued fractions with his unpublished 28 page report {\bf  \cite{1963-2}},  December, 1963. Continued fractions would play a role in his later research (cf., \cite{1972-2}). The main problem herein involves the finding of a lower bound on the real part of a function of infinitely many complex variables (though linear in each variable separately) when the variables are all assumed to lie in a given fixed disk in the complex plane. The lower bound thus obtained is used to show that the radii of univalence of some continued fractions are identical with their radii of starlikeness (see \cite{1963-2} for more details).

\subsection*{Asymptotic and oscillation theory of ordinary differential equations}\label{aot}\vskip0.25in

The article {\bf \cite{1951}} is Atkinson's first paper in the subject of (ordinary) differential equations: A two-term linear equation is considered with solutions and coefficients having values in a complex commutative Banach algebra. An approximate representation of the general solution is given to within a {\it little-oh} term.  In {\bf \cite{1951-2}} the author presents an asymptotic formula for the number of zeros of a solution (along with approximations to the solutions themselves) of a classical two-term smooth Sturm-Liouville equation on a semi-axis using a modified Pr\"{u}fer transformation (a technique that he used with ingenuity on many occasions). This is his first paper on classical Sturm-Liouville theory, a topic which he returned to many times in his career. In {\bf \cite{1954}} Atkinson studies the asymptotic behavior of solutions of linear perturbations of superlinear two-term ordinary differential equations of the second order along with some results on the oscillation and non-oscillation of its solutions. The paper {\bf \cite{1954-2}} is similar in spirit to \cite{1954} in that it deals in part with the asymptotic integration of  Emden-Fowler type equations. It also contains asymptotic results for linear equations extending results of Wintner (1949).

Atkinson's {\bf \cite{1955}} was a landmark in the study of oscillations of Emden-Fowler type differential equations and their extensions. He actually presents a necessary and sufficient condition for all solutions of said equation to oscillate on a half-axis provided the coefficient of the nonlinearity is positive. These results are still being extended to this day (see {\em e.g.,} \cite{Butler} and Section 2.3 in \cite{ABMbook}) for versions more general than Atkinson's. In {\bf \cite{1955-3}} he finds asymptotic estimates for the bounded solutions of oscillatory second order two-term asymptotically linear ordinary differential equations. The contribution {\bf \cite{1957-2}} is in the spirit of \cite{1954-2} in that a similar question is asked albeit under different hypotheses on the coefficients. Once again his interest is in the determination of solution asymptotics of a two-term second order linear ordinary differential equation under conditions somewhat different but akin to those in \cite{1957-2}.

Although {\bf \cite{1958}} is a technically difficult article to describe in few words the arguments there are related (according to Hartman) to similar ones used in a paper of Hartman and Wintner (1955). A real system of  fully nonlinear non-autonomous differential equations is considered on a semi-axis with a view at determining the uniform boundedness of its solutions on the given semi-axis from boundedness on compact subsets. It includes many results of the time. The write-up {\bf \cite{1968-1}} sees Atkinson attack the problem of asymptotic stability of periodic linear systems by reducing it to a case where a previous result of Levinson (1948) can be applied. His study of the area function in {\bf \cite{1972-5}} (akin to \cite{1960-2}) is used to find a Lyapunoff function for Emden-Fowler type equations thereby obtaining asymptotic results for its solutions in general cases. He returns once more to the Emden-Fowler equation in {\bf \cite{1978}}, a survey-type communication which announces a result reminiscent of his next paper \cite{1978-1}, but in the context of nonlinear equations.

There is little evidence anywhere that could have prepared researchers for the appearance and impact of {\bf \cite{1978-1}}. Even though I was likely the only other person to have read this work before its appearance in print (it was one of his sort-of-papers given to me as a student) he didn't give me any motivation for the reasons behind it, thus its existence goes beyond the realm of reason and established methods and since there is nothing like it anywhere (except possibly \cite{1954-2})  it most probably belongs to the world of magic. It appeared out of nowhere. In an attempt to gain greater insight into the asymptotics of  Emden-Fowler type equations he considers a special linear second order equation with a parameter appearing quadratically and shows that the set of all parameter values for which this equation admits a solution not tending to zero at infinity is an additive group. This is then used to derive a necessary and sufficient condition for certain nonlinear equations to have solutions tending to zero at a predetermined rate.

The study of the non-oscillatory behavior of solutions of non-homogeneous second order linear equations is the object of {\bf \cite{1980}}. For example the authors show that if the homogeneous equation is in the limit-circle case, then the corresponding non-homogeneous equation will have at most one non-oscillatory solution provided the square of the non-homogeneous term as an average value that grows at a given rate (in this regard, the techniques are reminiscent of Atkinson's article, \cite{1949-1}). 

The study of the stability of differential-difference equations is undertaken in {\bf \cite{1982-3}} using a modified version of Gronwall's inequality adapted to this situation. Later he and Haddock produce {\bf \cite{1983}}, a paper in which is found a discussion of the asymptotics of a first-order nonlinear delay differential equation on a half-line. Although his first published paper in differential equations with delays appeared in 1982, he certainly had an interest in this field at least as far back as 1976, when he discussed some problems in the area with the author.

During the mid-eighties there was a rush to prove the so-called {\it Maximum Eigenvalue Conjecture} attributed to W.T. Reid. The idea behind this conjecture was  to produce a vector systems analog of the classical Fite-Wintner theorem for the oscillation of second order real scalar differential equations using Wintner's notion of a ``conjugate point" as a zero for vector solution of a symmetric second order linear differential system, and by replacing the usual integral condition by the largest eigenvalue of the integral in question.  The publication {\bf \cite{1985}} was one of the many papers that contributed to this problem, eventually solved by M.K. Kwong {\it et al} shortly after.

His visit to the University of Alberta in the mid-eighties led to the joint-paper with Jack Macki {\bf \cite{1986-2}} where the authors return to an old problem dealing with zero asymptotics of solutions of second order two-term linear equations with a potential that is unbounded at infinity (see also \cite{1978-1}). His only report ever to appear in a Chinese journal is {\bf \cite{1987-4}}, where the authors obtain sufficient conditions for all the solutions of a linear first order delay differential equation to tend to zero at infinity. In {\bf \cite{1987-6}} he undertook the study of non-oscillatory Emden-Fowler equations in a book that we had edited, {\bf \cite{1987-5}} in the references, in commemoration of the sesquicentennial anniversary of the appearance of the pioneering paper by C.F. Sturm (1836) on boundary value problems for second order linear differential equations. More information is provided on the unique positive solution appearing in a previous article by Kaper and Kwong dealing with a free boundary problem.

\subsection*{Functional analysis and operator theory}\label{fa}\vskip0.25in
Derick's interests quickly shifted from the theory of the Riemann zeta-function to differential equations and functional analysis around 1948-1949. His first contribution to functional analysis {\bf \cite{1949-4}} includes the introduction of the definition of a symmetric operator on a general complex Banach space and its ramifications, {\em e.g.,} the existence of orthoprojectors and the orthogonality of eigenfunctions. The motivation for the next publication {\bf \cite{1951-1}} in this field is unclear, although it may have arisen as a result of some research in the former Soviet Union by Nikols'ski\u{i}. It was his first article in a language (Russian) other than his native English and it became rather important in the Soviet Union although its existence in the western world would surface many years later. There he introduces the notion of a {\it generalized Fredholm operator} and proves an index theorem. Then, using the theory of normed rings, he creates an abstract version of a theorem of Mihlin on singular integral operators. The applications are to singular integral equations of Mihlin type. The impact of \cite{1951-1} is discussed in Section~\ref{five}.

The next article in the subject, {\bf \cite{1952-1}}, is an excursion into the spectral theory of polynomial pencils. There he presents a general result to the effect that the set of values of $\lambda$ of a polynomial $T(\lambda)$ with compact coefficients such that $T(\lambda)f=0$ is an isolated set in the complex plane. This seems to have led the route to his future studies into multiparameter spectral theory and his extensions of Klein's oscillation theorem. Curiously enough, this article was immediately followed by B. Nagy's alternate proof (see \cite{nagy}) of the same result using complex function theory methods. Thus, it may be that Bela Sz.Nagy was the referee of Atkinson's \cite{1952-1}. In {\bf \cite{1952-2}} he extends a result of Yosida by weakening a condition on the {\it big-oh} term. The paper {\bf \cite{1953}} introduces the notion of a {\it relatively regular} operator along with in-depth study of the null spaces of such an operator and its adjoint. These results extend previous results of Halilov (1949), Gohberg (1952) and his own, \cite{1951-1} (see also Section~\ref{five}).

By 1963 his interests included Wiener-Hopf theory, perhaps motivated by the school of M.G. Kre\u{i}n and his classic paper of 1958 on the same subject. The first paper in his contribution to the evolution of this theory is {\bf \cite{1963-1}}. Within we find a thorough study of linear operators on an algebra into itself satisfying a generalized {\it integration by parts} formula. It is used to derive numerous applications in discrete and continuous Wiener-Hopf theory. In addition, a generalized factorization theorem is also presented, the theory being closely allied to Kre\u{i}n's own developments of the subject. 

After a long diversion into differential equations he returns to an old problem in {\bf \cite{1975}} dealing with his work \cite{1949-4} of 1949. Basically, much of the theory of the preceding paper is  extended here to {\it normaloid} operators on a Banach space (see the paper for detailed definitions).

\subsection*{Functions of a complex variable}\label{fcv}\vskip0.25in

Many of the results appearing in the subsection on the Riemann Zeta-function and Number Theory may well have been included here. We choose instead to include here those results that do not fit in the theory of numbers but deal with functions of a complex variable in the large. In this topic we find {\bf \cite{1961}}, an article that deals with finding a nice upper bound for the mini-max of a product of complex numbers thereby improving on some results of Erd\"{o}s and Szekeres (1961).

That same year, 1961, saw the publication of his first paper on sums of powers of complex numbers {\bf \cite{1961-1}}, a topic that he would return to on another occasion \cite{1969}. Incidentally, the first talk of the graduate lecture series {\it Conversations in Mathematics} held regularly at the University of Toronto from 1975-1979, was Derick Atkinson's who spoke exactly on this topic during the Fall of 1975. In \cite{1961-1} he answers a question posed by Tur\'{a}n in 1953. Atkinson proves that given $n$ points in the unit disk (the first of which is on the unit circle) arranged in order of decreasing moduli, the modulus of the sum of their $k^{th}$-powers is bounded below by an absolute constant (Derick gets $1/6$). In {\bf \cite{1969}} he returns to the problem of estimating the sums of powers of $n$ complex numbers previously considered in {\bf \cite{1961-1}} by improving the lower bound to $\pi/8$ and then actually refining this estimate with an asymptotic one, as $n \to \infty$.

\subsection*{Orthogonal polynomials}\label{op}\vskip0.25in

His insatiable interests in analysis led him, in the midst of his excursions into operator theory and functional analysis (see above), to a paper on orthogonal polynomials {\bf \cite{1952}}, written in German. There he finds a lower bound on the number of sign changes of a linear combination of orthogonal polynomials where the polynomial of lowest degree is subject to a constrained minimum condition in a weighted Lebesgue space $L^{p}_{\rho}(a,b)$, with $1 < p < \infty$. His article {\bf  \cite{1954-3}} is an investigation into systems of polynomials obtained by orthogonalization of a sequence of powers of $x$ with respect to a given positive measure. The resulting polynomials are compared to the classical ones and show similar characteristics. The next paper in this group {\bf \cite{1955-2}}, appeared in 1955. Here he obtains asymptotic estimates for a class of orthogonal polynomials whose derivatives vanish at $\pm 1$. In his paper {\bf \cite{1956}} Derick returns to a  problem first encountered in \cite{1952}. The basic difference here is that $p=\infty$ is allowed, the result being a separation theorem for the roots of two consecutive polynomials arranged in order of ascending degrees. The next article in this series deals with a theory, following Szeg\"{o} (1939), of orthogonal polynomials in several variables. He does this by introducing multiparameter notation and defining the polynomials as solutions of three-term-type recurrence relations from which one obtains the ``orthogonality" as a consequence. It is here that tensor products are first used to visualize the polynomials; this device will be crucial later on in his studies of multiparameter spectral theory (cf., \cite{1968}). 

In {\bf \cite{1972-2}} he returns to the study of orthogonal polynomials indirectly by asking for a polynomial $p_n(z)$ of degree n such that $p_n(z)F(z)$ has a Laurent expansion in which some given sequence of n coefficients must all vanish. Here $F$ is represented by a Stieltjes integral over a finite interval (reminiscent of a Stieltjes transform; cf., also \cite{1954-3}). Given two Tchebytcheff systems on a finite interval he shows in {\bf \cite{1973}} that there is a unique polynomial that gives a best ``bilateral" approximation. For such a polynomial he also proves a separation property of Sturm type for the zeros of two consecutive terms in the sequence as is done classically.  His next paper {\bf \cite{1981-1}}  in this area was with Everitt. In it they extend the classic result of Bochner (1929) on the classification of orthogonal polynomial solutions to real second order differential equations by allowing the coefficients to depend on the degree of the polynomial concerned using Stieltjes transforms.

\subsection*{Spectral theory of ordinary differential operators}\label{std}\vskip0.25in

His first paper in spectral theory in the strict sense is {\bf \cite{1957}}. Atkinson finds upper and lower bounds on the smallest positive eigenvalue of two-term Sturm-Liouville operators of even degree with a weight-function under some boundary conditions resembling periodic ones. The results include the cases where the weight function is non-negative and where it is periodic. The ``Bremmer series" paper {\bf \cite{1960-1}} consists in the finding of sufficient conditions for the expansion of the solutions of a smooth one-dimensional wave equation into terms $u_n(x)$ which can be interpreted as ``contributions generated after $n$ internal reflections in this medium" (Bremmer). 

The appearance of the monograph {\bf \cite{1964}} brought Atkinson international acclaim. It was so successful that it was translated into Russian shortly thereafter with the addition of an appendix by the noted M.G. Kre\u{i}n and I. S. Kac. It is full of ideas and new approaches to the field of boundary value problems (of the second order mainly)  for both ordinary differential and difference equations (in the guise of three-term recurrence relations). Results are presented for matrix systems as well and for Volterra-Stieltjes integral equations with a view at a unification of differential and difference equations (cf., \cite{ABMbook} for a somewhat dated survey of such an approach). For the reader who is unfamiliar with this text we strongly recommend Erdelyi's review in \cite{erde}. With a preliminary copy in hand, Richard Bellman foretold in 1962, that this book ``will become a classic" which of course, it did. 

The long term impact left by the publication of the paper {\bf \cite{1968}} cannot be underestimated. This was the first paper in a subject discussed earlier by Klein and others on the subject of multiparameter spectral theory, an article that led in due course to a renaissance in the area, and one which led to the publication of his  first volume on the subject \cite{1972}. He cast the subject in an abstract algebraic framework, one that would allow for generalizations to abstract operators with many parameters. Current work will likely see the publication of his long awaited second volume in the future (that is, the continuation of \cite{1972}). It seems that Atkinson already had the book {\bf \cite{1972}} in mind when he wrote \cite{1968}. In it he sets the framework/preliminaries for his second volume that would deal with applications to differential equations, however, this volume never appeared in print (although an incomplete draft appears to exist among his papers).

Atkinson was always interested in the Weyl-classification of second order ordinary differential equations and the associated spectral theory. His first paper in this area seems to be {\bf \cite{1972-1}} joint with Des. Evans. The authors give criteria for the classical two-term singular Sturm-Liouville operator to be of {\it limit-point}-type in a weighted Lebesgue space of square-integrable functions. When I first met Derick in 1974 and asked him about his research interests he jokingly considered himself a member of the {\it Limit-Point, Limit-Circle Classification Society} given his all-consuming interest in this area, at that time. Thus, in {\bf \cite{1972-3}} he improves on some results of Everitt and Giertz. Using the area/volume preserving device (see \cite{1960-2} and \cite{1972-5})  he obtains limit-point criteria for nonlinear second order equations including those of Emden-Fowler type in {\bf \cite{1973-2}}, a class of equations that he returned to constantly in his studies of asymptotic and oscillation properties of differential equations.

The extension of Weyl's limit-point, limit-circle theory for singular second order boundary problems with an indefinite weight-function (but with a left-definite left-hand side)  was the subject of his joint paper {\bf \cite{1974-1}} with Everitt and his fourth doctoral student, K. S. Ong. It became the source of much future investigation in particular, in papers by H. Langer and co-authors (1970's). The solo paper {\bf \cite{1974-2}} became an instant classic in that he considered not only second order but actually even order singular differential expressions under very general assumptions on their coefficients. Using criteria on sequences of intervals rather than on a half-line as was usually done he obtains new ``limit-n" criteria for these equations. Most of the results of the time are embedded in this majestic paper which won him the Makdougall-Brisbane Prize of the Royal Society of Edinburgh for the best paper covering the two-year period, 1974-1976. As the limit-point, limit-circle classification society  gathered steam around the world Derick, having tackled the even order theory in \cite{1974-2}, now decided to look at a corresponding theory for {\it squares} of differential operators restricting himself to the second order case and its square, in {\bf \cite{1976-1}}. 

The preprint {\bf \cite{1976}} is a preliminary version of the joint paper {\bf \cite{1977}} with Eastham and McLeod where the authors tackle the problem of classifying the $L^2$-solutions of a two-term Sturm Liouville operator on a half-axis with a potential of the form $q(x) = x^{\alpha}p(x^{\beta})$, where $p$ is periodic and all exponents are positive. In true Atkinson style he uses his knowledge of the field in order to further the study of the limit-point, limit-circle theory of Sturm-Liouville equations with {\it several} parameters thus opening a new direction for future researchers in {\bf \cite{1977-1}}. In a joint paper with Everitt {\bf \cite{1978-2}}, they discuss the optimal value of the parameter such that a classical singular Sturm-Liouville operator on a half-line will have a solution that is square-integrable at infinity. In a subsequent paper {\bf \cite{1980-1}}, he gives estimates on those values of the parameter for which the solutions of either a scalar or vector Sturm-Liouville equation with a weight on a half-line decay (or grow) exponentially fast and the distance of these values from the essential spectrum of  the corresponding operator,  a topic he returned to later with Des Evans in {\bf \cite{1982}} for complex valued coefficients in a more general framework. His work {\bf \cite{1980}} on the location of Weyl circles, whose existence are guaranteed by the {\it Weyl classification} (of a second order singular equation into limit-point, limit circle)  consists in an improvement of the asymptotics of the m-function based on a previous paper by Everitt and Halvorsen (1978).

Although the 1970's were the key period for Derick's research into the limit-point, limit-circle theory he still pursued this interest well into the 1980's with {\bf \cite{1981-3}}, a paper that dealt with extensions of the theory to matrix systems, difference equations and possibly complex-valued coefficients. According to sources close to Derick he wrote this paper on the evening before his talk at the conference (with very little, if any, sleep). Once the Weyl classification more or less satisfactorily complete so far as Derick was concerned, he turned to a study of the Titchmarsh-Weyl m-coefficient during much of the 1980's with a handful of papers, {\bf \cite{1982-1}, \cite{1984}, \cite{1988}, \cite{1988-4}, \cite{1988-5}} and the paper with Charles Fulton {\bf \cite{1999}}, published after his illness. The following articles with Fulton namely, {\bf \cite{1982-2}, \cite{1983-1}}, and {\bf \cite{1984-1}}, all dealt with limit-circle problems and eigenvalue asymptotics.

Immersed within the most prolific period of his career, the decade of the 1980's, we find a mere handful of papers in his favorite subject, the spectral theory of ordinary differential operators. Motivated in part by a paper of the author in 1982, he and an undergraduate student, David Jabon, produced numerical evidence for the existence of complex eigenvalues for indefinite (normally termed {\it non-definite} these days) Sturm-Liouville problems on a finite interval by studying (both numerically and theoretically) the case where the potential is constant and the weight function is piecewise constant and changes sign once. This was followed later (see {\bf \cite{1987}}) by a rigorous proof of the J\"{o}rgens Conjecture, dating from 1964, on the distribution of the zeros (and so of the eigenvalues) of generally weighted Sturm-Liouville problems on a finite interval, which included the non-definite case in particular. 

Motivated by trendsetting work in mathematical oceanography by Boyd in the early eigthies and John Adam in magnetohydrodynamics, Atkinson and others proceed to establish a rigorous foundation for real Sturm-Liouville problems with an interior pole in {\bf \cite{1988-1}} and {\bf \cite{1988-2}} (see also \cite{1988-5}). His last sole authored publication was {\bf \cite{1991-3}} where he finds asymptotic estimates for radially symmetric solutions of a semilinear Laplacian in the unit ball in three-space and obtains asymptotic formulae for the eigenvalues where the corresponding eigenfunctions satisfy an oscillation theorem of sturm-type on a finite open interval.

Although unable to oversee the final draft of the papers {\bf \cite{1993-0}} and {\bf \cite{1994}} because of his disability, the idea behind it was to define and study a Titchmarsh-Weyl m-function for a second order linear differential equation in which the coefficients are allowed to vary with the eigenvalue parameter. This complements some work in magnetohydrodynamics by John Adam as mentioned earlier. These ideas were later generalized to two-dimensional block operator matrices in {\bf \cite{1994}}.

\subsection*{Partial differential equations}\label{pdes}\vskip0.25in

His first contribution to differential equations was in the area of {\it partial} differential equations: The article {\bf \cite{1949-3}} includes a proof of the uniqueness of the solution of the classical wave equation outside a surface $\mathcal{S}$ given its boundary values on $\mathcal{S}$ provided the Sommerfeld radiation condition (along with boundedness)  is satisfied at infinity. After an excursion of many years into ordinary differential equations he met up with Bert Peletier, who would become his long time collaborator, with whom he wrote {\bf \cite{1971-1}}. It deals with the existence of a unique weak nonnegative solution of the problem of a laminar flow for a gas flowing through a homogeneous porous material. Most of the papers of his later years (1986-1992) would deal with questions arising from partial differential equations many of which could be reduced to ordinary differential equations, Atkinson's playground.

In his {\bf \cite{1974}} paper with Peletier he looks at the problem of finding solutions to a boundary value problem corresponding to a nonlinear diffusion equation with concentration-dependent diffusion coefficient using the method of successive approximations after reduction to the one variable case. The study of positive solutions asymptotic to zero of a semilinear Laplacian in connection with the Emden-Fowler equation is the subject of the article {\bf \cite{1986}}, an investigation that was completed in {\bf \cite{1986-4}} for the two-dimensional case. In this connection the paper {\bf \cite{1986-3}} deals with a perturbation of the classic Emden-Fowler equation with a critical exponent and the subsequent development of the theory of its positive solutions as the perturbation tends to zero. The asymptotic behavior of radial solutions of the special equation in the title of {\bf \cite{1986-1}} is considered in the light of the exponent $p$  not being necessarily rational. In {\bf \cite{1987-1}} the authors perturb the exponent in a semilinear Laplacian and find the asymptotic behavior of the resulting positive radial solution as the exponent tends to zero. In {\bf \cite{1987-3}} Atkinson and Peletier provide ordinary-differential-equations proofs of classic results on the existence of positive solutions of nonlinear elliptic equations considered earlier by Br\'{e}zis and Nirenberg along with some extensions (this line of research was very {\it hot} at the time).

The contents of the article {\bf \cite{1988-6}} include results on the non-existence of radial solutions that change sign in the unit ball for a nonlinear Laplacian with a parameter appearing linearly in $R^N$ where $N=4,5,6$. The subject of the existence of positive radial solutions to nonlinear elliptic equations is reconsidered in {\bf \cite{1988-7}} in the light of special eigenvalue-eigenfunction asymptotics for the Dirichlet problem on the unit ball. The subject of the existence of positive solutions to a nonlinear (also in its principal part) second order ordinary differential equation that are bounded away from zero at infinity is treated in {\bf \cite{1988-8}}.  However, in  {\bf \cite{1988-9}} the authors initiate a new kind of investigation into the existence of positive radial solutions (that are not necessarily functions) of the prescribed mean curvature equation in n-space. The contribution {\bf \cite{1988-11}} describes the (positive) solution set of a semilinar Laplacian with a nonlinearity that has critical growth at infinity. The existence of oscillatory solutions to the m-Laplacian in the critical exponent case is considered in {\bf \cite{1988-12}}. 

The article {\bf \cite{1989-2}} and its predecessor {\bf \cite{1989-1}}, consider the problem of the existence of finite vertical tangent lines of radial solutions of the prescribed mean curvature equation by reduction to an ordinary equation (see also \cite{1988-9} and the later {\bf \cite{1991-1}}). {\bf \cite{1989-3}} is considered a preprint by Peletier (2004) and the results therein appeared elsewhere (likely in \cite{1988-7}). {\bf \cite{1989-4}} is the preprint of the refereed publication {\bf \cite{1990-1}} in which the authors seek to address the problem of the oscillatory solutions of the eigenvalues and eigenfunctions of the Dirichlet problem associated with a semilinear Laplacian in n-space. Once again, the main tool is reduction to an ordinary equation. In {\bf \cite{1990}} the authors consider the existence and non-existence of oscillatory solutions  (or nodal solutions), the asymptotic estimate for their zeros and ultimately the asymptotics for the eigenvalues of elliptic problems with critical Sobolev exponents.

His eternal inspiration from the Emden-Fowler equation led him to write {\bf \cite{1991}} where an ordinary equation of Emden-Fowler type is considered with an eye at obtaining asymptotics of specific solutions (those initiating at zero with value $\gamma$ whose derivative is zero) as $\gamma$ tends to infinity. 

\subsection*{Miscellaneous topics}\label{mt}\vskip0.25in

The first article that is off the beaten track is likely {\bf \cite{1960}} written with two-co-authors, G. A. Watterson and P.A.P. Moran. Incidentally, Pat Moran was instrumental in bringing Atkinson to Canberra from Ibadan and the year of submission was also the year that Atkinson left for Toronto. The article deals with relationships between the sum of the entries of a non-negative rectangular matrix $A$ and the sum of the elements of $AA^TA$. An {\it area principle} is considered in {\bf \cite{1960-2}} for a two-dimensional time-dependent Hamiltonian system from which one can draw information of solution asymptotics for small values of a parameter. This area principle was reconsidered later in \cite{1972-5}. One question that Atkinson returned to a number of times under different covers consisted in the solvability of the equation $fx=y$ for $x$. In {\bf \cite{1961-2}} he considers the problem of the solvability of this equation where $f$ is a differentiable function on $n-$dimensional Euclidean space into itself. As a consequence he obtains an $n-$dimensional extension of the fundamental theorem of algebra. It is difficult to determine the rationale for his first paper in mathematical statistics {\bf \cite{1964-1}}, joint with two other authors. It is a basic contribution to decision making under complete ignorance that involves the existence of a minimax (likely Atkinson's contribution, in part). There are few researchers who appreciated inequalities more than Atkinson, - he created them and used them continually in his papers and on occasion published them for their own sake. One such example is the article {\bf \cite{1971}}. Using a simple argument he finds a {\it lower} bound  on the integral of a product of two (generally convex) functions in terms of the product of the individual integrals.

His love of inequalities (undoubtedly following the tradition of Hardy-Titchmarsh) produced the remarkable paper {\bf \cite{1972-4}} (cf., also \cite{1960},\cite{1961},\cite{1961-1},\cite{1971}). In this last paper dealing explicitly with inequalities he determines conditions on a general nonlinearity appearing in a two-term second order ordinary equation on a half-line under which solutions will be positive, unique, with some {\it a priori} bounds. The only published paper of a historical nature is {\bf \cite{1981-2}} dealing with Christoffel's work on shock waves where he discusses two key papers by Christoffel (1877). In the joint paper {\bf \cite{1988-10}} with Haddock, they discuss the well-posedness of systems of functional differential equations with possibly infinite delays and some applications. His two papers {\bf \cite{1997-0}} and {\bf \cite{milloux}} with the late \'{A}rpad Elbert reveal an interest into his extending various ideas from linear differential equations theory to the {\it half-linear} case well-known to have been a specialty of the Bihari-Elbert school.

Although numerical computation and/or programming seem at first to be incongruous with Atkinson's research interests, this can only be due in part, to his sincere modesty.  The legacy of works he left behind is riddled with BASIC or even Fortran code, some of the work dating back to early times, {\it ca.} mid-1960's. Thus, it should not come as a surprise that he would have a hand in contributing to numerical procedures as he did in the papers {\bf \cite{1989-5}} and {\bf \cite{leslie}} including a specific mention for his help in {\bf \cite{anes}}, a paper in anesthesia.

\section{Publications}\label{pubs}

We present a list of the publications of  F. V. Atkinson arranged in chronological order. Articles published in a common year are usually ordered by date of submission or receipt (when available) not by date of appearance. This gives a better perspective on the fields that he was working on at the time.

\section{Impact of the works of F.V. Atkinson}\label{five}

We cite some references to works of F.V. Atkinson. We cannot make any claim as to completeness; it is meant to complement the text by providing information as to the impact of some of the publications referred to above in various areas of mathematics and where appropriate or even possible, to give an update to work in that area.\vskip0.25in

His highly influential work on the Riemann zeta-function (see above) still has an impact on the area. For recent work referencing Atkinson's papers from the period (1939-1950) see \cite{bala1}, \cite{bbc},\cite{chak1},\cite{jutila1},\cite{gotze},\cite{hb1},\cite{hb2}, \cite{ivic1} among many others. His paper \cite{1939-1} on the partition function was cited in the survey paper \cite{gupta}. His papers on Tur\'{a}n's problem (cf., \cite{1961-1},\cite{1969}) still gather interest, see \cite{cheer}. One of his results, \cite{1951-2}, was even used in connection with black-hole geometries, \cite{black}! At the time of writing, distinct papers referencing at least one of Atkinson's papers number over 1200.\vskip0.25in

Curiously enough, one of the most enduring papers of his early days is actually \cite{1949-3}, see the thorough survey paper \cite{sommerfeld}.  The main result in this paper led to what is now called the {\it Atkinson-Wilcox Theorem}, partly in his honor. Its applications are widespread, for example in ellipsoidal geometry \cite{fokas}, in thermoelasticity \cite{dassios2}, the theory of acoustics \cite{farfield}, etc. Besides his seminal work on the Sommerfeld radiation condition his next most influential mathematical work was in the field of the theory of operators (with \cite{1951-1}) where his generalization of Fredholm operators and the corresponding index theory led to what is now known as {\it Atkinson's Theorem}, see \cite{kato}, \cite{halmos}. Some of  his other classic work in operator theory (see \cite{1949-4}, \cite{1975})  was recently extended by Zhuravlev in \cite{zhur}. More recent work \cite{1994}, has seen applications to medicine \cite{colitis}, \cite{colitis2}; the paper with his son Leslie, \cite{leslie}, in psychology \cite{gaps},\cite{honjo}, his earlier contribution to medicine \cite{anes}, and other areas mostly dealing with mathematics.\vskip0.25in

Apart from his published articles his two books \cite{1964} and \cite{1972} have had a profound influence on the field, the first by treating discrete and continuous boundary problems separately and then together along with their extensions, while the second by providing the algebraic formalism for multiparameter spectral theory in a long awaited sequel that unfortunately never appeared during his life. \vskip0.25in

\section{Acknowledgments}\label{ak}
I am deeply indebted to Dusja Atkinson, her son Leslie and her daughter Vivienne, for fascinating and illuminating conversations and for access to the vast literature left behind by F.V. Atkinson. 

My sincere thanks also go to Kathi and Philip Leah, friends of the Atkinson family; Prof. P.G. Rooney (University of Toronto); Michael Riordan (Archivist, St. Johns and The Queen's Colleges Archives at Oxford; Prof. W. Norrie Everitt (the University of Birmingham, U.K.); Prof. Don B. Hinton (University of Tennessee); the American Mathematical Society and the editors of  Mathematical Reviews; Prof. Bert Peletier (University of Leiden, The Netherlands); Prof. Reinhard Mennicken (University of Regensburg, Germany); Alex Scott and John Gallehawk (Bletchley Park archivists); Sam Howison, old time friend of Derick Atkinson; Prof. Lazlo Hatvani (University of Szeged, Hungary); Prof. Chandler Davis (University of Toronto); Simon Bailey (Keeper of the Archives, University of Oxford, Bodleian Library).


\begin{thebibliography}{199}
\bibitem{1939} {\it The mean-value of the zeta-function on the critical line}\\
Quarterly J. Math. Oxford Series, {\bf 10} (38)  (1939), 122-128.
\bibitem{1939-1} {\it A summation formula for $p(n)$, the partition function}\\
J. London Math. Soc. {\bf 14} (1939), 175-184. \quad{(MR 1-40b)}
\bibitem{1941}{\it The mean-value of the zeta-function on the critical line}\\
Proc. London Math. Soc.(2), {\bf 47} (1941) 174-200.\quad{(MR 3-70a)}
\bibitem{1941-1} {\it A divisor problem}\\
Quarterly J. Math. Oxford Series, {\bf 12} (48)  (1941), 193-200.\quad{(MR 3-269d)}
\bibitem{1948} {\it The Abel summation of certain Dirichlet series}\\
Quarterly J. Math. Oxford Series, {\bf 19} (73)  (1948), 59-64.\quad{(MR 9-508a)}
\bibitem{1948-1} {\it \"{U}ber die Stirlingsche Reihe}\\
Comment. Math. Helvetici, {\bf 21} (4) (1948), 332-335.\quad{(MR 10-32c)}
\bibitem{1948-2} {\it A mean-value property of the Riemann zeta-function}\\
J. London Math. Soc. {\bf 23} (1948), 128-135.\quad{(MR 10-182d)}
\bibitem{1949-1} {\it The mean-value of the Riemann zeta-function}\\
Acta Math. {\bf 81} (1949), 353-376.\quad{(MR 11-234d)}
\bibitem{1949-2} (with Lord Cherwell) {\it The mean-values of arithmetical functions}\\
Quarterly J. Math. Oxford Series, {\bf 20} (1949), 65-79.\quad{(MR 11-15b)}
\bibitem{1949-3} {\it On Sommerfeld's ``Radiation Condition"}\\
Philos. Magazine,\, Ser. 7, {\bf 40} (1949), 645-651.\quad{(MR 10-714a)}
\bibitem{1949-4} {\it Symmetric Linear Operators on a Banach Space }\\
Monatshefte f\"{u}r Math. {\bf 53} (4) (1949), 278-297. \quad{(MR 11-525e)}
\bibitem{1950} {\it The Riemann zeta function}\\
Duke Math. J. {\bf 17} (1) (1950), 63-68.\quad{(MR 11-162b)}
\bibitem{1951} {\it Asymptotic properties of a differential equation}\\
Actas Acad. Sci. Lima {\bf 14} (1951), 28-33.\quad{(MR 13-653a)}
\bibitem{1951-1} {\it The normal solubility of linear equations in normed spaces\,\,}(in Russian)\\
Mat. Sbornik, NS, {\bf 28} (70) (1951), 3-14.\quad{(MR 13-46d)}
\bibitem{1951-2} {\it On second order linear oscillators}\\
Revista (Serie A) Mat. y Fisica Teorica, Universidad Nacional de Tucum\'{a}n, {\bf 8} (1-2) (1951), 71-87. \quad{(MR 14-50a)}
\bibitem{1952} {\it \"{U}ber die Nullstellen gewisser extremaler Polynome}\\
Archiv der Math. {\bf 3} (2) (1952), 83-86. \quad{(MR 14-269b)}
\bibitem{1952-1} {\it A spectral problem for completely continuous operators}\\
Acta Math. Acad. Sci. Hungaricae, {\bf 3} (1-2) (1952), 53-60. \quad{(MR 14-478h)}
\bibitem{1952-2} {\it On a theorem of Yosida}\\
Proc. Japan Acad. {\bf 28} (1952), 327-329.\quad{(MR 14-381c)}
\bibitem{1953} {\it On relatively regular operators}\\
Acta Sci. Math. (Szeged), {\bf 15} (8) (1953), 38-56.\quad{(MR 15-134e)}
\bibitem{1954} {\it On linear perturbation of non-linear differential equations}\\
Can. J. Math. {\bf 6} (4) (1954), 561-571.\quad{(MR 16-701e)}
\bibitem{1954-2} {\it The asymptotic solution of second-order differential equations}\\
Ann. Mat. Pur. App. (4) {\bf 37} (1954), 347-378.\quad{(MR 16-701f)}
\bibitem{1954-3} {\it On lacunary and other orthogonal polynomials}\\
Revista (Serie A) Mat. y Fisica Teorica, Universidad Nacional de Tucum\'{a}n, {\bf 10} (1-2) (1954), 95-110. \quad{(MR 17-32b)}
\bibitem{1955} {\it On second-order non-linear oscillations}\\
Pacific J. Math. {\bf 5} (1) (1955), 643-647.\quad{(MR 17-264e)}
\bibitem{1955-2} {\it On orthogonal polynomials with extrema at the ends of the orthogonality interval}\\
Monats. f. Math. {\bf 59} (4) (1955), 323-330.\quad{(MR 17-607e)}
\bibitem{1955-3} {\it On asymptotically linear second-order oscillations}\\
J. Rational Mech. Analysis {\bf 4} (5) (1955), 769-793.\quad{(MR 17-264d)}
\bibitem{1956} {\it On polynomials with least weighted maximum}\\
Proc. Amer. Math. Soc. {\bf 7} (2) (1956), 267-270.\quad{(MR 18-126a)}
\bibitem{1957} {\it Estimation of an eigenvalue occuring in a stability problem}\\
Math. Zeitschrift {\bf 68} (1957), 82-99.\quad{(MR 19-1052a)}
\bibitem{1957-2} {\it Asymptotic formulae for linear oscillations}\\
Proc. Glasgow Math. Assoc. {\bf 3} (3) (1957), 105-111.\quad{(MR 22-9655)}
\bibitem{1958} {\it On stability and asymptotic equilibrium}\\
Annals Math. (2), {\bf 68} (3) (1958), 690-708.\quad{(MR 20-7137)}
\bibitem{1960} (with G.A. Watterson and P.A.P. Moran) {\it A matrix inequality}\\
Quarterly J. Math. Oxford Second Series, {\bf 11} (42)  (1960), 137-140.\quad{(MR 22-9502)}
\bibitem{1960-1} {\it Wave propagation and the Bremmer series}\\
J. Math. Analysis Appl. {\bf 1} (1960), 255-276.\quad{(MR 23-B1846)}
\bibitem{1960-2}{\it A constant area principle for steady oscillations}\\
J. Math. Analysis Appl. {\bf 1} (2) (1960), 133-144.\quad{(MR 22-6912)}
\bibitem{1961} {\it On a problem of Erd\"{o}s and Szekeres}\\
Canad. Math. Bull. {\bf 4} (1) (1961), 7-12\quad{(MR 23-A3722)}
\bibitem{1961-1} {\it On sums of powers of complex numbers}\\
Acta Math. Acad. Sci. Hungaricae {\bf 12} (1-2) (1961), 185-188.\quad{(MR 23-A3714)}
\bibitem{1961-2} {\it The reversibility of a differentiable mapping}\\
Canad. Math. Bull. {\bf 4} (2) (1961), 161-181.\quad{(MR 23-A3219)}
\bibitem{1963} {\it Boundary problems leading to orthogonal polynomials in several variables}\\
Bull. Amer. Math. Soc. {\bf 69} (3) (1963), 345-351.\quad{(MR 26-3950)}
\bibitem{1963-1} {\it Some aspects of Baxter's functional equation}\\
J. Math. Analysis Appl. {\bf 7} (1) (1963), 1-30.\quad{(MR 27-5135)}
\bibitem{1963-2} {\it A value-region problem occuring in the theory of continued fractions}\\
MRC Technical Summary Report \# 419, Univ. Wisconsin Math. Res. Center, December 1963, 28 pp.
\bibitem{1964} {\it Discrete and Continuous Boundary Problems}\ (Book)\\
Academic Press, New York, 1964. \quad{(MR 31-416)}\\
(Also translated into Russian by Izdat ``MIR", Moscow, 1968, with an Addendum by I. S. Kac and M. G. Kre\u{i}n, {\it R-functions: analytic functions mapping upper half-plane to itself.}\quad{(MR 39-4473)}
\bibitem{1964-1}  (with J. D. Church and B. Harris)\ {\it Decision procedures for finite decision problems under complete ignorance}\\
Ann. Math. Stat. {\bf 35} (4) (1964), 1644-1655.\quad{(MR 29-6593)}
\bibitem{1968} {\it Multiparameter spectral theory}\\
Bull. Amer. Math. Soc. {\bf 74} (1)  (1968), 1-27.\quad{(MR 36-3145)}
\bibitem{1968-1} {\it On asymptotically periodic linear systems}\\
J. Math. Analysis Appl. {\bf 24} (3) (1968), 646-653.\quad{(MR 39-3095)}
\bibitem{1969} {\it Some further estimates concerning sums of powers of complex numbers}\\
Acta Math. Acad. Sci. Hungaricae {\bf 20} (1-2) (1969), 193-210.\quad{(MR 39-443)}
\bibitem{1971} {\it An inequality}\\
Publ. Fac. D'\'{E}lectrotechnique Univ. Belgrade (S\'{e}rie Math. Phys.) (357-380) (1971), 5-6.\\ \quad{(MR 45-6998)}
\bibitem{1971-1} (with L.A. Peletier) {\it Similarity profiles of flows through porous media}\\
Arch. Rat. Mech. Analysis {\bf 42} (5) (1971), 369-379.\quad{(MR 48-12984)}
\bibitem{1972} {\it Multiparameter Eigenvalue Problems, Vol 1- Matrices and Compact Operators} (Book)\\
Academic Press, New York, (1972), xi, 209 p.\quad{(MR 56-9291)}
\bibitem{1972-1} (with W. Desmond Evans) {\it On solutions of a differential equation which are not of integrable square}\\
Math. Zeitschrift {\bf 127} (4) (1972), 323-332.\quad{(MR 48-4390)}
\bibitem{1972-2}  {\it Orthogonal polynomials and lacunary approximants}\\
in {\it Orthogonal Polynomials and their Applications}, P. Graves-Morris, G.A. Baker Jr, eds.  Proc. Conf. Kent, 17-21 July, 1972;\\
Acad. Press, New York, London, (1972), 75-79.
\bibitem{1972-3} {\it On some results of Everitt and Giertz}\\
Proc. Roy. Soc. Edinb., {\bf 71A} (13) (1972/73), 151-158.\quad{(MR 48-4391)}
\bibitem{1972-4} {\it On second-order differential inequalities}\\
Proc. Roy. Soc. Edinb., {\bf 72A} (8) (1972/73), 109-127.\quad{(MR 53-5996)}
\bibitem{1972-5} {\it The area function for non-linear second-order oscillations}\\
Proc. Roy. Soc. Edinb., {\bf 72A} (10) (1972/73), 135-147.\quad{(MR 54-13208)}
\bibitem{1973} {\it Bilateral approximations and \v{C}eby\v{c}ev systems}\\
Archiv der Math. {\bf 24} (3) (1973), 297-302.\quad{(MR 47-9138)}
\bibitem{1973-2} {\it Nonlinear extensions of limit-point criteria}\\
Math. Zeitschrift {\bf 130} (4) (1973), 297-312.\quad{(MR 47-8963)}
\bibitem{1974} (with L. A. Peletier) {\it Similarity solutions of the nonlinear diffusion equation}\\
Arch. Rat. Mech. Analysis {\bf 54} (4) (1974), 373-392.\quad{(MR 49-9298)}
\bibitem{1974-1} (with W.N. Everitt and K.S. Ong) {\it On the $m$-coefficient of Weyl for a differential equation with an indefinite weight function}\\
Proc. London Math. Soc., (3rd Series) {\bf 29} (1974), 368-384.\quad{(MR 53-8546)}
\bibitem{1974-2} {\it Limit-$n$ criteria of integral type}\\
Proc. Roy. Soc. Edinb., {\bf 73A} (11) (1974/75), 167-198.\quad{(MR 52-6090)}
\bibitem{1975} {\it The polynomial-normaloid property for Banach-space operators}\\
Monatshefte Math. {\bf  79} (1975), 273-283.\quad{(MR 51-11167)}
\bibitem{1976-1} {\it Asymptotic integration and the $L^2$-classification of squares of second-order differential operators}\\
Quaestiones Math. {\bf 1} (1976), 155-180.\quad{(MR 56-12389)}
\bibitem{1977} (with M.S.P. Eastham and J. B. McLeod) {\it The limit-point, limit-circle nature of rapidly oscillating potentials}\\
Proc. Roy. Soc. Edinb., {\bf 76A} (1977), 183-196.
\bibitem{1977-1} {\it Deficiency-index theory in the multi-parameter Sturm-Liouville case}\\
in {\it Differential Equations}, Proceedings from the Uppsala 1977 International Conference on Differential Equations, Symposium \# 7, Uppsala (1977), 1-10.\quad{(MR 57-16789)}
\bibitem{1978} {\it Asymptotics of certain non-linear oscillatory ordinary differential equations}\\
in {\it Proceedings of the VIIIth. International Conference on Nonlinear Oscillations}, Prague (1978), 117-121.
\bibitem{1978-1} {\it A stability problem with algebraic aspects}\\
Proc. Roy. Soc. Edinb., {\bf 78A} (1978), 299-314.\quad{(MR 58-11634)}
\bibitem{1978-2} (with W.N. Everitt) {\it Bounds for the point spectrum for a Sturm-Liouville equation}\\
Proc. Roy. Soc. Edinb., {\bf 80A} (1978), 57-66.\quad{(MR 80g:34021)}
\bibitem{1980} (with R.C. Grimmer and W.T. Patula) {\it Nonoscillatory solutions of forced second-order linear equations, II}\\
Annali Mat. Pura ed Appl. (4), {\bf 76} (1980), 299-317.\quad{(MR 82i:34032)}
\bibitem{1980-1} {\it Exponential behaviour of eigenfunctions and gaps in the essential spectrum}\\
in {\it Lecture in Mathematics 827}, Proc. 1978 Dundee International Conference on Differential Equations \\
Springer-Verlag, New York, (1980), 1-24.\quad{(MR 82f:34022)}
\bibitem{1981} {\it On the location of Weyl circles}\\
Proc. Roy. Soc. Edinb., {\bf 88A} (1981), 345-356.
\bibitem{1981-1} (with W.N. Everitt) {\it Systems of orthogonal polynomials satisfying second-order differential equations}\\
in {\it E.B. Christoffel}, P.L. Butzer and F. Feh\'{e}r eds., Proc. International Christoffel Symposium, Aachen (1979) \\
Birkh\"{a}user Verlag, Basel (1981), 173-181.\quad{(MR 83j:33008)}
\bibitem{1981-2} {\it Christoffel's work on shock waves}\\
in {\it E.B. Christoffel}, P.L. Butzer and F. Feh\'{e}r eds., Proc. International Christoffel Symposium, Aachen (1979) \\
Birkh\"{a}user Verlag, Basel (1981), 718-720.\quad{(MR 84c:01027)}
\bibitem{1981-3} {\it A class of limit-point criteria}\\
in {\it Spectral Theory of Differential Operators}, I.W. Knowles and R.T. Lewis Eds., North-Holland Publ. \# 55 (1981), 13-35.\quad{(MR 83c:34023)}
\bibitem{1982} (with W.D. Evans) {\it On the exponential behaviour of eigenfunctions and the essential spectrum of differential operators}\\
Proc. Roy. Soc. Edinb., {\bf 92A} (1982), 271-300.\quad{(MR 84a:34027)}
\bibitem{1982-1} {\it On the asymptotic behaviour of the Titchmarsh-Weyl $m$-coefficient and the spectral function for scalar second-order differential expressions}\\
in {\it Lecture in Mathematics 964}, Proc. 1982 Dundee International Conference on Differential Equations \\
Springer-Verlag, New York, (1982), 1-27.\quad{(MR 85b:34022)}
\bibitem{1982-2} (with C.T. Fulton) {\it Some limit-circle eigenvalue problems and asymptotic formulas for eigenvalues}\\
in {\it Lecture in Mathematics 964}, Proc. 1982 Dundee International Conference on Differential Equations \\
Springer-Verlag, New York, (1982), 28-55.\quad{(MR 85e:34016)}
\bibitem{1982-3} (with J.R. Haddock and O.J. Staffans) {\it Integral inequalities and exponential convergence of solutions of differential delay equations with boundary delay}\\
in {\it Lecture in Mathematics 964}, Proc. 1982 Dundee International Conference on Differential Equations \\
Springer-Verlag, New York, (1982), 56-68.\quad{(MR 84g:34123)}
\bibitem{1983} (with J. R. Haddock) {\it Criteria for asymptotic constancy of solutions of functional differential equations}\\
J. Math. Analysis Appl.,  {\bf 91} (2) (1983), 410-423.\quad{(MR 84f:34096)}
\bibitem{1983-1} (with C. T. Fulton) {\it Asymptotic formulae for eigenvalues of limit circle problems on a half-line}\\
Annali Mat. Pura Appl., (IV) {\bf 135} (1983), 363-398.\quad{(MR 85j:34107)}
\bibitem{1984} {\it On bounds for Titchmarsh-Weyl $m$-coefficients and for spectral functions for second-order differential operators}\\
Proc. Roy. Soc. Edinb., {\bf 97A} (1984), 1-7.\quad{(MR 85j:35032)}
\bibitem{1984-1} (with C. T. Fulton) {\it Asymptotics of Sturm-Liouville eigenvalues for problems on a finite interval with one limit-circle singularity, I}\\
Proc. Roy. Soc. Edinb., {\bf 99A} (1984), 51-70.\quad{(MR 86f:34044)}
\bibitem{1984-2} (with D. Jabon) {\it Indefinite Sturm-Liouville problems}\\
in {\it Proc. 1984 Workshop on the Spectral Theory of Sturm-Liouville Differential Operators}, Argonne National Laboratory, ANL-84-73 (1984), 31-45.
\bibitem{1985} (with Hans G. Kaper and Man Kam Kwong) {\it An oscillation criterion for linear second-order differential systems}\\
Proc. Amer. Math. Soc., {\bf 94} (1) (1985), 91-96.\quad{(MR 86h:34028)}
\bibitem{1986} (with L. A. Peletier) {\it Ground states of $-\triangle u = f(u)$ and the Emden-Fowler equation}\\
Arch. Rat. Mech. Analysis {\bf 93} (2) (1986), 103-127.
\bibitem{1986-1} (with Lambertus-Adrianus Peletier) {\it Sur les solutions radiales de l'\'{e}quation $\triangle u +(1/2)x \cdot \bigtriangledown u + (1/2)\lambda u + |u|^{p-1}u =0.$}\\
C. R. Acad. Sci. Paris (S\'{e}rie I) {\bf 302} (3) (1986), 99-101.\quad{(MR 87j:35124)}
\bibitem{1986-2} (with J.W. Macki) {\it On regular growth and asymptotic stability}\\
Rocky Mountain J. Math. {\bf 16} (1986), 111-117.\quad{(MR 87m:34043)}
\bibitem{1986-3} (with L. A. Peletier) {\it Emden-Fowler equations involving critical exponents}\\
Nonlinear Analysis, Th. Meth. Appl., {\bf 10} (8) (1986), 755-776.\quad{(MR 87j:34039)}
\bibitem{1986-4} (with L. A. Peletier) {\it Ground states and Dirichlet problems for $-\Delta u = f(u)$ in $R^2$}\\
Arch. Rat. Mech. Analysis {\bf 96} (1986), 147-165.\quad{(MR 87k:35080)}
\bibitem{1987} (with A. B. Mingarelli) {\it Asymptotics of the number of zeros and of the eigenvalues of general weighted Sturm-Liouville problems}\\
J. Reine Ang. Math., {\bf 375-376} (1987), 380-393.\quad{(MR 88d:34023)}
\bibitem{1987-1} (with L.A. Peletier) {\it Elliptic equations with nearly critical growth}\\
J. Diff. Eq., {\bf 70} (3) (1987), 349-365.\quad{(MR 89e:35054)}
\bibitem{1987-3} (with L. A. Peletier) {\it Elliptic equations with critical exponents}\\
in {\it Nonlinear Parabolic equations: Qualitative Properties of Solutions}, Pitman Res. Notes in Math. 149, L. Boccardo and A. Tesei, eds. Longman and Wiley, (1987), 13-21.\quad{(MR 90d:35095)}
\bibitem{1987-4} (with Zhang Shunian) {\it Asymptotic behavior of solutions of linear delay equations}\\
Acta Math. Sinica, n.s. Vol. 3, (4) (1987), 289-300.\quad{(MR 89c:34083)}
\bibitem{1987-5} (with W. F. Langford and A. B. Mingarelli eds.)\,\, {\it Oscillation, Bifurcation and Chaos} (Book)\\
CMS Conf. Proc. Volume 8, AMS  Providence, R.I., (1987) xv, 740 p.
\bibitem{1987-6} (with Chen Shaozhu) {\it On the nonoscillation of an Emden-Fowler equation}\\
in {\it Oscillation, Bifurcation and Chaos}, (F.V. Atkinson, W. F. Langford and A. B. Mingarelli eds.) Canadian Mathematical Society Proceedings Vol. 8, AMS    Providence, R.I., (1988) 43-55.\quad{(MR 88k:34033)}
\bibitem{1988} {\it On the order of magnitude of Titchmarsh-Weyl functions}\\
Diff. and Integral Equations, {\bf 1} (1)  (1988), 79-96.
\bibitem{1988-1}  (with W. N. Everitt and A. Zettl) \,\, {\it Regularization of a Sturm-Liouville problem with an interior singularity using quasi-derivatives}\\
Diff. and Integral Equations, {\bf 1} (1988), 213-221.
\bibitem{1988-2}  {\it Asymptotics of an eigenvalue problem involving an interior singularity}\\
in {\it Proc. of the Focused Research Program on Spectral Theory and Boundary Value Problems: vol. 2, Singular Differential Equations}, Argonne National Laboratory, ANL-87-26, (1988), 1-18.
\bibitem{1988-4}  {\it Estimation of the Titchmarsh-Weyl function $m(\lambda)$ in a case with an oscillating leading coefficient}\\
in {\it Proc. of the Focused Research Program on Spectral Theory and Boundary Value Problems: vol. 2, Singular Differential Equations}, Argonne National Laboratory, ANL-87-26, (1988), 19-43.
\bibitem{1988-5}  (with C.T. Fulton) {\it Asymptotics of the Titchmarsh-Weyl $m$-Coefficient for non integrable potentials}\\
in {\it Proc. of the Focused Research Program on Spectral Theory and Boundary Value Problems: vol. 2, Singular Differential Equations}, Argonne National Laboratory, ANL-87-26, (1988), 79-103.
\bibitem{1988-6} (with H. Br\'{e}zis and L.A. Peletier) {\it Solutions d'\'{e}quations elliptiques avec exposant de Sobolev critique qui changent de signe}\\
C.R. Acad. Sci. Paris, S\'{e}rie I, {\bf 306} (1988), 711-714.\quad{(MR 89k:35088)}
\bibitem{1988-7} (with L.A. Peletier) {\it Large solutions of elliptic equations involving critical exponents}\\ Asymptotic Analysis {\bf 1} (1988), 139-160.\quad{(MR 89i:35004)}
\bibitem{1988-8} (with L.A. Peletier) {\it On non-existence of ground states}\\
Quart. J. Math. Oxford Series, (2), {\bf 39} (1988), 1-20.\quad{(MR 89b:34007)}
\bibitem{1988-9}  (with L.A. Peletier and J. Serrin) {\it Ground states for the prescribed mean curvature equation: the supercritical case}\\
in {\it Nonlinear Diffusion Equations and their Equilibrium States, I}, W.M. Ni, L.A. Peletier and J. Serrin eds., Mathematical Sciences Research Institute Publ. \#12, Springer-Verlag, New York, (1988), 51-74.\quad{(MR 89j:35051 )}
\bibitem{1988-10} (with J.R. Haddock) {\it On determining phase spaces for functional differential equations}\\
Funkcialaj Ekvacioj {\bf 31} (1988), 331-347.\quad{(MR 90f:45009)}
\bibitem{1988-12} {\it Asymptotics of eigenfunctions for some nonlinear elliptic problems}\\
in {\it Proc. International Conf. on Theory and Applications of Differential Equations, Columbus, Ohio, 1988}, A.R. Aftabizadeh ed. Ohio University Press, Vol.1 (1989), 26-49.\quad{(MR 91k:35086)}
\bibitem{1989} {\it Self-intersecting solutions of the prescribed mean curvature equation}\\ in {\it Proc. of the Focused Research Program on Spectral Theory and Boundary Value  Problems: Nonlinear Differential Equations, Vol. 4}, Argonne National Laboratory ANL-87-26 (August 1989), pp. 1-20.
\bibitem{1989-2} (with L.A. Peletier) {\it Bounds for vertical points of solutions of prescribed mean curvature type equations I}\\ 
Proc. Roy. Soc. Edinb., {\bf 112A} (1989), 15-32.\quad{(MR 90m:35059)}
\bibitem{1989-5}  (with A. M. Krall, G. K. Leaf and A. Zettl) {\it On the numerical computation of eigenvalues of Sturm-Liouville problems with matrix coefficients}\\
in {\it Proc. of the Focused Research Program on Spectral Theory and Boundary Value Problems: Vol. 3, Linear Differential Equations and Systems}, Argonne National Laboratory, ANL-87-26, (1989), 21-38.
\bibitem{leslie} (with L. Atkinson, B. Quarrington and J. J. Cyr)\, {\it Differential classification in school refusal.}\\
British Journal of Psychiatry, {\bf 155} (1989 ), 191-195.
\bibitem{1990} (with H. Br\'{e}zis and L.A. Peletier)  {\it Nodal solutions of elliptic equations with critical Sobolev exponents}\\
J. Differential Equations {\bf 85} (1) (1990), 151-170.\quad{(MR 91e:35021)}
\bibitem{1990-1} (with L.A, Peletier) {\it Oscillations of solutions of perturbed autonomous equations, with an application to nonlinear elliptic eigenvalue probems involving critical Sobolev exponents}\\
Differential and Integral Equations {\bf 3} (1990), 401-433.\quad{(MR 91c:34075)}
\bibitem{1991} {\it Asymptotic integration and a borderline case in the theory of elliptic equations involving critical exponents}\\
in {\it Proc. International Conference on Differential Equations, Dundee, July 1990}, Pitman Res. Notes Math. Ser. 254, Longman Sci. Tech., Harlow (1991), 1-24.\quad{(MR 93e:34076)}
\bibitem{1991-2} (with L.A. Peletier and J. Serrin) {\it Estimates for vertical points of solutions of prescribed mean curvature equations, II}\\ 
Asymptotic Analysis {\bf 5} (1992), 283-310.\quad{(MR 93c:35034)}
\bibitem{1991-3} {\it Higher approximations to eigenvalues for a nonlinear elliptic problem}\\
in {\it Nonlinear diffusion equations and their equilibrium states}, Proc. 1989 Conf. Gregynod, {\it Prog. Nonlinear Differential Equations Appl. \# 7} \\
Birkhauser, Boston (1991), 39-69.\quad{(MR 93f:35174)}

\section{Miscellaneous articles}\label{mas}
This section contains a list miscellaneous papers, reports, papers that have appeared elsewhere, and manuscripts authored by F.V. Atkinson (and/or co-authors), sorted in chronological form. The list also contains references to handwritten manuscripts for which there is no traceable date of origin and papers published during Derick Atkinson's long illness (1992-2002). We note that Derick was unable to contribute to any research during this period so the tacit assumption is made that research arising from joint authorship originated prior to 22 August, 1992.
\bibitem{thesis} {\it The analytic theory of numbers}\\
Dissertation for the D. Phil, Oxford University, Trinity Term 1939, iii, 63 p.
\bibitem{nodate-1} {\it Asymptotic properties of a differential equation}\, (Offprint)\\
Acad. Nacional Ciencias Exactas, Fisicas y Naturales, 1-6  (Post 1949).
\bibitem{1956-0} {\it Frontiers of Pure Mathematics}\\
Inaugural Lecture, Canberra University College, Canberra, Australia, 19th July, 1956. 17pp.
\bibitem{1976} (with M. Eastham and J. McLeod) \,\, {\it The limit-point, limit-circle nature of rapidly oscillating potentials}\\
MRC Technical Summary Report \# 1676, Univ. Wisconsin Math. Res. Center, September 1976, 30 pp.
\bibitem{1981-0} (with S.N. Zhang) {\it Asymptotic behavior and structure of solutions for equation $x^{\prime}(t) = p(t)(x(t) - x(t-1)).$} \\
Unpublished.
\bibitem{1987-2} (with L. A. Peletier) {\it Radial similarity solutions of a parabolic equation} (preprint)\\
in {\it Nonlinear Parabolic equations: Qualitative Properties of Solutions}, Pitman Res. Notes in Math. 149, L. Boccardo and A. Tesei, eds. Longman and Wiley, (1987), 5-12.\quad{(MR 90d:35134)}
\bibitem{1988-11} (with L.A. Peletier) {\it Elliptic equations with critical growth when $N \geq 3$ and $N=2$} (preprint) \\ Pitman Res. Notes in Math. 181, Longman and Wiley, (1988), 1-12.\quad{(MR 90i:35101)}
\bibitem{1989-1} (with L.A. Peletier) {\it Bounds for vertical points of solutions of prescribed mean curvature type equations} (preprint)\\ in {\it Proc. of the Focused Research Program on Spectral Theory and Boundary Value  Problems: Nonlinear Differential Equations, Vol. 4}, Argonne National Laboratory  ANL-87-26 (August 1989), pp. 21-42.
\bibitem{1989-3} (with L.A. Peletier) {\it Large solutions of elliptic equations involving critical exponents}(preprint)\\ in {\it Proc. of the Focused Research Program on Spectral Theory and Boundary Value  Problems: Nonlinear Differential Equations, Vol. 4}, Argonne National Laboratory ANL-87-26 (August 1989), pp. 43-72.
\bibitem{1989-4} (with L. A. Pelletier)\,\, {\it Oscillation of solutions of perturbed autonomous equations, with an application to nonlinear elliptic eigenvalue problems involving critical Sobolev exponents} (preprint)\\
Argonne National Laboratory MCS-P91-0789, (1989).
\bibitem{1991-1} (with L.A. Peletier and J. Serrin) {\it Bounds for vertical points of solutions of prescribed mean curvature equations II} (preprint)\\ 
Mathematical Institute, University of Leiden, Report W91-02, January 1991. 
\bibitem{1993-0}  (with H. Langer, R. Mennicken) {\it Sturm-Liouville problems with coefficients which depend analytically on the eigenvalue parameter}\\
Acta Sci. Math. (Szeged) {\bf 57} (1993), 25 - 44. \quad{(MR 94m:34060)}
\bibitem{1994}  (with H. Langer, R. Mennicken, A.A. Shkalikov) {\it The essential spectrum of some matrix operators}\\
Math. Nach. {\bf 167} (1994), 5 - 20. \quad{(MR 95f:47007)}
\bibitem{1995} (with H.G. Kaper and M.K. Kwong) {\it Asymptotics for a free boundary problem}\\
Methods and Applications of Analysis {\bf 2} (4) (1995), 466-474.\quad{(MR 97a:35017)}
\bibitem{1997-0} (with \'{A}. Elbert) {\it Extension of Prodi-Trevisan theorem to a half-linear differential equation}\, 
(1997), Preprint.
\bibitem{1999}  (with C.T. Fulton) {\it Asymptotics of the Titchmarsh-Weyl $m$-coefficient for non-integrable potentials, IV}\\
Proc. Roy. Soc. Edinb. {\bf 129A} (1999), 663-683.\quad{(MR 2000f:34047)}
\bibitem{milloux} (with \'{A}. Elbert) {\it An extension of Milloux's theorem to half-linear differential equations}\\
in {\it Proceedings of the 6th. Colloquium on the Qualitative Theory of Differential Equations}, Electron. J. Qual. Theory Differ. Equ.(Szeged), (8) (2000), 10p.


\bibitem{bala1} R. Balasubramanian and K. Ramachandra, \,\, {\it Mean-value of the Riemann Zeta-Function on the critical line}\\
Proc. Indian Acad. Sci.-Math. Sci. {\bf 93}  (1984), 101-107.
\bibitem{bbc} B. C. Berndt, \,\, {\it Identities involving coefficients of a class of Dirichlet series. 4}\\
Trans. Amer. Math. Soc. {\bf 149} (1) (1970), 179.
\bibitem{Butler} G. J. Butler, \,\, {\it On the oscillatory behavior of a second order nonlinear differential equation}\\
Ann. Mat. Pura ed Appl. Ser. 4, {\bf 105} (1975), 73-91.
\bibitem{chak1} I.C. Chakravarty, \,\, {\it Certain properties of a pair of secondary zeta-functions}\\
J. Math. Analysis Appl. {\bf 35} (3) (1971), 484.
\bibitem{cheer}  A. Y. Cheer and D. A. Goldston, {\it  Tur\'{a}n's Pure Power Sum Problem}\\
Math. Comp. {\bf  65} (1996), 1349-1358.
\bibitem{erde} A. Erdelyi, in {\it Mathematical Reviews}, {\bf 31}, No. 416 (MR 0176141)
\bibitem{dassios2} F. Cakoni and G. Dassios, \,\, {\it The Atkinson-Wilcox theorem in thermoelastcity}\\
Quart. Applied Math. {\bf 57} (4) (1999), 771-795.
\bibitem{fokas} G. Dassios, \,\, {\it The Atkinson-Wilcox theorem in ellipsoidal geometry}\\
J. Math. Analysis Appl. {\bf 274} (2)  (2002), 828-845.
\bibitem{colitis2} Delaney CP, Remzi FH, Gramlich T, et al.,\,\,{\it Equivalent function, quality of life and pouch survival rates after ileal pouch-anal anastomosis for indeterminate and ulcerative colitis}\\
Annals of Surgery, {\bf  236} (1) (2002), 43-48.
\bibitem{jutila1} M. Jutila, \,\, {\it Riemann's Zeta-Function and the divisor problem}\\
Ark. Math. {\bf 21} (1) (1983), 75-96.
\bibitem{colitis} Geboes K, De Hertogh G, \,\, {\it Indeterminate colitis}\\
Inflammatory Bowel Diseases, {\bf 9 }(5) (2003), 324-331.
\bibitem{gotze} F. Gotze, \,\, {\it Mean properties of Riemann zeta function}\\
J. Math. Soc. Japan {\bf 19} (4) (1967), 426.
\bibitem{gupta} H. Gupta, \,\, {\it Partitions - A survey}\\
J. Res. National Bureau of Standards-B, Math. Sci. {\bf 74} (1) (1970), 1. 
\bibitem{halmos} P.R. Halmos, \,\, {\it A Hilbert Space Problem Book}\\
Van Nostrand Reinhold, Princeton NJ (1967).
\bibitem{hatvani} L. Hatvani, {\it On stability properties of solutions of second order differential equations}\\
in {\it Proceedings of the 6th. Colloquium on the Qualitative Theory of Differential Equations}, Electron. J. Qual. Theory Differ. Equ.(Szeged), (11) (2000), 6p.
\bibitem{hb1} D. R. Heath-Brown, \,\, {\it Fourth power moment of the Riemann zeta-function}\\
Proc. London Math. Soc. {\bf 38} (1979), 385-422.
\bibitem{hb2} D. R. Heath-Brown, \,\, {\it Distribution and moments of the error term in the Dirichlet divisor problem}\\
Acta Arithmetica {\bf 60} (4) (1992), 389-415.
\bibitem{honjo} Honjo S, Nishide T, Niwa S, et al., \,\, {\it School refusal and depression with school inattendance in children and adolescents: Comparative assessment between the Children's Depression Inventory and somatic complaints}\\
Psychiat. Clin. Neuros. {\bf 55} (6) (2001), 629-634
\bibitem{ivic1} A. Ivic, \,\, {\it On the error term for the fourth moment of the Riemann zeta-function}\\
J. London Math. Soc. (2) {\bf 60} (1) (1999), 21-32.
\bibitem{kato} T. Kato, \,\, {\it Perturbation Theory of Linear Operators}\\
Springer-Verlag; Reprint Ed. edition (February 1995) , 619 pp.
\bibitem{gaps}Kearney CA, \,\, {\it Bridging the gap among professionals who address youths with school absenteeism: Overview and suggestions for consensus}\\
Professional Psychology- Research and Practice, {\bf 34} (1) (2003),  57-65.
\bibitem{black} J. Louko and S. Makela, \,\,{\it Area spectrum of the Schwarzschild black hole}\\
Phys. Rev. D, {\bf 54} (8) (1996), 4982-4996.
\bibitem{ABMbook}A. B. Mingarelli, \,\, {\it Volterra-Stieltjes Integral Equations and Generalized Ordinary Differential Expressions}\\
Lecture Notes in Mathematics 989, Springer Verlag, New York, 1983, xiv, 317 p.
\bibitem{jwp} J. W. Porter, \,\, {\it Generalized Titchmarsh-Linnik divisor problem}\\
Proc. London Math. Soc. (24) (1972), 15.
\bibitem{farfield} C. Prabavathi and C.P. Vendhan, \,\, {\it Determination of far-field pattern of rigid scatterers using independent finite element method and eigenfunction expansion.2. Non-axisymmetric scattering}\\
J. Vibration and Acoustics, Trans. ASME, {\bf 118} (4) (1996), 583-590.
\bibitem{sommerfeld} Schott S.H., \,\, {\it Eighty years of the Sommerfeld radiation condition}\\
Historia Math. {\bf 19} (4) (1992), 385-401.
\bibitem{nagy}  B. Sz. Nagy, {\it On a spectral problem of Atkinson}\\
Acta Math. Acad. Sci. Hungaricae, {\bf 3} (1-2) (1952), 61-66.
\bibitem{tit1} E.C. Titchmarsh, {\it The Theory of the Riemann Zeta Function}, Oxford University Press, $1951$.
\bibitem{anes} A.H.C. Walker and R.J. Stout, {\it The effects of anaesthesia upon fallopian tubal mobility}\\
J. Obstetrics and Gynaecology of the British Empire, {\bf 59} (1) (1952), 1-24.
\bibitem{weinstein1} L. Weinstein, \,\, {\it Zeros the Artin L-series of a cubic field on the critical line}\\
J. Number Theory {\bf 11} (2) (1979), 279-284.
\bibitem{zhur} V.F. Zhuravlev, \,\, {\it Generalization of a theorem of FV Atkinson for a case of one class of normally solvable operators in Banach spaces}\\
Dokl. Akad. Nauk. {\bf 358} (2) (1998), 157-159.

\section {Miscellaneous references}\label{mrefs}
\bibitem{abstract} {\it Abstracts of Dissertations for the degree of Doctor of Philosophy}\\
University of Oxford Committee for Advanced Studies Vol. 12, Oxford at the Clarendon Press, 1940, v, 353 p.
\bibitem{bautin}N.N. Bautin, {\it On the number of limit cycles which appear with the variation of coefficients from an equilibrium position of focus or center type}\\
Mat. Sbornik (N.S.) (30) (72) (1952) 181-196,\\
Translated into English by F.V. Atkinson-published by the AMS in 1954 as Translation Number 100, 396-413; reprint, AMS Transl.(1) 5 (1962), 396-413.
\bibitem{cms} {\it CMS-SMC} $1945-1995$, Vol.1, {\it Mathematics in Canada}, Canadian Mathematical Society, Ottawa, $1995$, viii, 406 p.
\bibitem{earl} The Earl of Birkenhead,\,\,{\it The Prof in Two Worlds, - The official life of Professor F. A. Lindemann, Viscount Cherwell}\\ Collins, St. James's Place, London, 1961, 383 p.
\bibitem{genea} http://genealogy.math.ndsu.nodak.edu/
\bibitem{hardy} G.H. Hardy, {\it A Mathematician's Apology}, Cambridge University Press, reprinted 1992.
\bibitem{lord} Frederick Alexander Lindemann (alias Lord Cherwell), 1886-1957, scientist, became the personal assistant to Sir Winston Churchill during the war years and also head of the Statistical section in his Admiralty. He received the honorific title of Baron in 1941 and Viscount in 1957 likely due in part to his service to England (e.g., in helping to establish the UK Atomic Energy Authority, his appointment as paymaster general of the UK, co-creator of the {\it Nernst-Lindemann theory of specific heats}, etc).
\bibitem{oliphant} Sir Marcus 'Mark' Laurence Elwin Oliphant (1901-2000), physicist, was instrumental during the Manhattan Project that helped to build the atomic bomb in the USA. Being one of Ernest Rutherford's (Nobel Laureate in Chemistry, 1908)  most distinguished students he lived in Australia most of his life.
\bibitem{archives} Michael Riordan, Archivist, St.Johns and the Queen's Colleges, Oxford University; {\it personal correspondence}, 2003-2004.
\bibitem{pgr} P.G. Rooney, {\it Obituary: Frederick Valentine Atkinson}\\
Austral. Math. Soc. Gaz. {\bf 30} (2) (2003), 100.\quad{(MR 1982044)}
\bibitem{bp} Alex Scott and John Gallehawk, {\it personal correspondence}, Bletchley Park, March 2003.
\bibitem{selberg} A. Selberg, \,\,{\it Harmonic analysis and discontinuous groups in weakly symmetric Riemannian spaces with applications to Dirichlet series}\\
J. Indian Math. Soc. (n.S.) {\bf 20} (1956), 47-87.
\bibitem{cpsnow} C. P. Snow, \,\,{\it Variety of Men}\\
Macmillan \& Co., London 1967, 224 pp.
\end{thebibliography}
\end{document}